\documentclass[leqno,final]{siamltex}
\usepackage{latexsym,amssymb,amsfonts, amsmath}
\usepackage{color}
\usepackage{enumerate}
\usepackage{mathrsfs}

\DeclareMathOperator*{\argmin}{argmin}
\DeclareMathOperator*{\argmax}{argmax}
\DeclareMathOperator*{\sign}{sign}
\newtheorem{remark}{\bf Remark}

\renewcommand{\b}[1]{\mbox{\boldmath $#1$}}
\title{Asymptotic stability of the Wonham filter: ergodic and non-ergodic signals
}

\author{Peter Baxendale
\thanks{
Department of Mathematics, University of Southern California, Los
Angeles, CA 90089-1113, USA ({\tt baxendal@math.usc.edu}). The
research of this author was supported by  ONR Grant
N00014-96-1-0413.} \and Pavel Chigansky
\thanks{
Dept. Electrical Engineering-Systems, Tel Aviv University, 69978
Tel Aviv, Israel ({\tt pavelm@eng.tau.ac.il}). } \and Robert
Liptser
\thanks{
Dept. Electrical Engineering-Systems, Tel Aviv University, 69978
Tel Aviv, Israel ({\tt liptser@eng.tau.ac.il}).} }

\begin{document}
\maketitle
% ----------------------------------------------------------------
\begin{abstract}
Stability problem of the Wonham filter with respect to initial
conditions is addressed.
The case of ergodic signals is revisited
in  view of a gap in the classic work  of H. Kunita (1971). We
give new bounds for the exponential stability rates, which do not
depend on the observations. In the non-ergodic case, the stability
is implied by identifiability conditions, formulated explicitly in
terms of the transition intensities matrix and the
observation
structure.
\end{abstract}



\pagestyle{myheadings} \thispagestyle{plain} \markboth{P.
BAXENDALE, P. CHIGANSKY, R. LIPTSER}{STABILITY OF THE WONHAM
FILTER}

% ----------------------------------------------------------------
\section{Introduction}\label{sec-1}

The optimal filtering estimate of a signal from the record of
noisy observations is usually generated by a nonlinear recursive
equation subject to the signal a priori distribution. If the
latter is unknown and the filtering equation is initialized  by an
arbitrary initial distribution, the obtained estimate is
suboptimal in general. From an applications point of view, it is
important to know whether such estimate becomes close to the
optimal one at least after enough time elapses. This property of
filters to forget the initial conditions is far from being obvious
and in fact generally remains an open and challenging problem.

In this paper, we consider the filtering
setting for signals
 with a finite state space.
Specifically, let $X=(X_t)_{t\ge 0}$ be a continuous time
homogeneous Markov chain observed via
\begin{eqnarray}\label{1.1}
Y_t = \int_0^t h(X_s)ds + \sigma W_t
\end{eqnarray}
with the Wiener process $W=(W_t)_{t\ge 0}$, independent of $X$,
some bounded function $h$, and $\sigma\neq0$.

We assume that $X_t$ takes values in the
finite alphabet
$\mathbb{S}= \{a_1,...,a_n\}$ and admits several ergodic classes.
Namely,
$$
\mathbb{S}=\Big\{\underbrace{a^1_1,\ldots,a^1_{n_1}}_{\mathbb{S}_1},
\ldots,\underbrace{a^m_1,\ldots,a^m_{n_m}}_{\mathbb{S}_m}\Big\},
$$
where the subalphabets $\mathbb{S}_1,\ldots,\mathbb{S}_m$ are noncommunicating
in the sense that for any $i \neq j$ and $t\ge s$
\begin{equation}\label{1.2'}
P\big(X_t\in\mathbb{S}_j|X_s\in\mathbb{S}_i\big)=0.
\end{equation}
So, unless $m=1$, $X_t$ is a compound Markov chain with the
transition intensities matrix
\begin{equation}\label{LA}
\Lambda=
  \begin{pmatrix}
    \Lambda_1 & 0 & 0
\\
    0 & \Lambda_2 & 0
\\
\ldots & \ldots & \ldots
\\
0 & 0 & \Lambda_m
  \end{pmatrix}
\end{equation}
of $m$ ergodic classes and is not ergodic itself.

The filtering problem consists in computation  of the conditional
distribution,
\[
\pi^\nu_t(1)=P(X^\nu_t=a_1|\mathscr{Y}^\nu_{[0,t]}),\ldots,
\pi^\nu_t(n)=P(X^\nu_t=a_n|\mathscr{Y}^\nu_{[0,t]}),
\]
where $\mathscr{Y}^\nu_{[0,t]}$ is the filtration, generated by
$\{Y^\nu_s,0\le s\le t\}$ satisfying the usual conditions
(henceforth, the superscript $\nu$ is
used to emphasize that the distribution of $X_0$ is $\nu$).

The  vector-valued random process $\pi^\nu_t$ with
entries
$\pi^\nu_t(1),\ldots,\pi^\nu_t(n)$ is generated by the Wonham
filter \cite{W} (see also \cite[Chap. 9]{LSI})
\begin{equation}\label{1.3n}
\begin{split}
\pi^\nu_0= \ &\nu,
\\
d\pi^\nu_t= \ &\Lambda^*\pi^\nu_t dt +
\sigma^{-2}\big(\diag(\pi_t^\nu)-\pi_t^\nu(\pi_t^\nu)^*\big)
h(dY^\nu_t-h^*\pi^\nu_tdt),
\end{split}
\end{equation}
where ${\diag}(x)$ is the
scalar matrix with the diagonal $x\in
\mathbb{R}^n$, $ h $ is the column vector with entries
$h(a_1),\ldots,h(a_n),$ and $^*$ is the transposition symbol. If
$\nu$ is unknown and some other distribution $\beta$ (on
$\mathbb{S}$) is used to initialize the filter, the ``wrong''
conditional distribution $\pi^{\beta\nu}_t$ is obtained:
\begin{equation}\label{1.4n}
\begin{split}
\pi^{\beta\nu}_0= \ &\beta,
\\
d\pi^{\beta\nu}_t = \ &\Lambda^*\pi^{\beta\nu}_t dt + \sigma^{-2}
\big(\diag(\pi_t^{\beta\nu})-\pi_t^{\beta\nu}(\pi_t^{\beta\nu})^*\big)
h (dY^\nu_t- h ^*\pi^{\beta\nu}_tdt).
\end{split}
\end{equation}

According to the
intuitive notion of stability, given at the
beginning of this section, the filter defined in \eqref{1.4n} is
said to be asymptotically stable if
\begin{equation}\label{sense'}
\lim_{t\to\infty}E\|\pi^\nu_t-\pi^{\beta\nu}_t\|=0,
\end{equation}
where $\|\cdot\|$ is the total variation norm.

If the state space of the Markov chain $X$ consists of  one
ergodic class ($m=1$), our setting is in the framework studied by
Ocone and Pardoux \cite{OP}. In this case, there exists the
unique invariant distribution $\mu$, so that
\begin{equation}\label{1.7m}
\lim\limits_{t\to\infty}\|S_t\gamma-\mu\|=0,
\end{equation}
where  $S_t$ is the semigroup corresponding to $X$ and $\gamma$ is
an arbitrary probability distribution on $\mathbb{S}$. Moreover
\begin{eqnarray}
\label{Kcond} \lim\limits_{t\to\infty}
\int_\mathbb{S}|\mathcal{S}_tf(x)-\mu(f)|d\mu(x)=0
\end{eqnarray}
holds for any bounded $f:\mathbb{S}\mapsto \mathbb{R}$. So, it may
seem that it remains only to assume
\begin{equation}\label{nullbeta}
\nu\ll\beta
\end{equation}
and allude to \cite{OP}. However, the proof of \eqref{sense'}
given in \cite{OP} uses as its central argument the uniqueness
theorem for the stationary measure of the filtering process
$\pi_t^\nu$ which appeared in the work of H. Kunita \cite{K}.
Unfortunately, the proof of this theorem  (Theorem 3.3  in
\cite{K}) contains a serious gap, as elaborated in the next
section.

Different approach to the stability analysis of the filters for ergodic
signals was initiated by Delyon and Zeitouni \cite{DZ}. The authors
studied  the top Lyapunov exponent of the filtering equation
$$
\gamma_\sigma(\beta',\beta'') = \varlimsup_{t\to
\infty}\frac{1}{t}\log\big\|\pi^{\beta'\nu}_t-\pi^{\beta''\nu}_t\big\|,
\quad \beta' \text{\ and\ } \beta''\text{\ distributions on\
}\mathbb{S},
$$
and show that $ \gamma_\sigma(\beta',\beta'')< 0 $
when $\Lambda$ and $h$ satisfy
certain conditions. Moreover the filter is found to be stable in
the low signal-to-noise regime: $
\varlimsup\limits_{\sigma\to\infty}\gamma_\sigma(\beta',\beta'')
\le \Re\big[\lambda^{\max}\big(\Lambda\big)] $ with
$\lambda^{\max}(\Lambda)$ being  the eigenvalue of $\Lambda$ with
the largest nonzero real part.

These results were further extended by Atar and Zeitouni
\cite{AZ1}, where it is shown that uniformly in $\sigma>0$ and $h$
\begin{eqnarray}\label{atar}
\gamma_{\sigma}(\beta',\beta'')\le -2\min_{p\ne
q}\sqrt{\lambda_{pq}\lambda_{qp}}, \ \text{a.s.},
\end{eqnarray}
and the
high signal-to-noise asymptotics are  obtained:
\begin{align*}
& \varlimsup_{\sigma\to 0}\sigma^2\gamma_\sigma \le
-\frac{1}{2}\sum_{i=1}^d \mu_i\min_{j\ne
i}\big[h(a_i)-h(a_j)\big]^2
\\
& \varliminf_{\sigma\to 0}\sigma^2\gamma_\sigma \ge
-\frac{1}{2}\sum_{i=1}^d \mu_i
\sum_{j=1}^d\big[h(a_i)-h(a_j)\big]^2
\end{align*}
where $\mu$ is the ergodic measure of $X$.

The method in \cite{AZ1} (and its full development in \cite{AZ2})
does not rely on \cite{K} and is based on the analysis of the Zakai
equation, corresponding to \eqref{1.3n} (see \eqref{zakai} below). The
analysis is carried out by means of the Hilbert projective metric
and the Birkhoff inequality, etc.; see
section \ref{sec-5} for more
details. This approach proved out its efficiency in several filtering
scenarios (see \cite{A}, \cite{BO2}, \cite{BKu-98}).

Other results and methods related to the filtering stability can
be found in \cite{AVZ}, \cite{BO}, \cite{BKu-99}, \cite{BKu-00},
\cite{Cerou}, \cite{DFM}, \cite{DMG}, \cite{DMG01}, \cite{COC},
\cite{LGM}, \cite{DMM01}, \cite{LGO01}, \cite{LGO02},  \cite{O1},
\cite{O2}. The linear Kalman$-$Bucy case, being the most understood,
is extensively treated by several authors: \cite{BK}, \cite{M},
\cite{MS}, \cite{DZ}, \cite{OP}, \cite{LBR}, \cite{LSII} (sections
14.6 and 16.2).

In the present paper, we consider both ergodic and non-ergodic
signals. Applying the technique from Atar and Zeitouni, \cite{AZ2}, we
show that in the ergodic case the asymptotic stability holds true
without any additional assumptions. In other words, the conclusion of H. Kunita \cite{K}
is valid in the specific case under consideration.

In view of the counterexample given in section \ref{sec-3}, it is
clear that in general $\gamma_\sigma$
may vanish at $\sigma=0$. So, it is interesting to
find out which ergodic properties of the signal are inherited by
the filter regardless of the specific observation structure. In
this connection we prove the inequality
$$
\varlimsup_{t\to\infty}\frac{1}{t}\log\|\pi^{\beta\nu}_t-\pi^\nu_t\|
\le -\sum_{r=1}^n\mu_r\min_{i\ne r}\lambda_{ri}.
$$
Since $\mu$ is the positive measure on
$\mathbb{S}$, unlike
\eqref{atar}, this bound remains negative if at least one row of
$\Lambda$ has all nonzero entries.

Also we give the nonasymptotic bound (compare with
\eqref{atar})
$$
\|\pi^\nu_t-\pi^{\beta\nu}_t\|\le C\exp\Big(-2 t\min_{p\ne
q}\sqrt{\lambda_{pq} \lambda_{qp}}\Big)
$$
with some positive constant $C$ depending on $\nu$ and $\beta$
only.

For the discrete time case,
related results can be found in Del Moral and
Guionnet \cite{DMG01} and Le Gland and Mevel \cite{LGM}. For example, in \cite{LGM}
the positiveness assumption for all
transition probabilities  is relaxed under certain constraints on
the observation process noise density.

In the case of nonergodic signal, $m>1$, we show that the
filtering stability holds true if the ergodic classes can be
identified via observations and the filter matched to each class
is stable. We formulate explicit sufficient identifiability
conditions in terms of $\Lambda$ and $h$.

The paper is organized as follows. In section \ref{sec-2}, we
introduce the necessary notations and clarify the role of condition
$\nu\ll\beta$ in the filtering stability (Proposition
\ref{pro-1.1}). This section also gives a link to the gap in
Kunita's proof \cite{K}, while in section \ref{sec-3} the filtering
setting is described for which the stability fails and the
gap becomes evident.

The main results are formulated in section
\ref{sec-4} and proved in sections \ref{sec-5} and \ref{sec-6}.

\section{Preliminaries and connection to the gap in \cite{K}}\label{sec-2}

\subsection{Notations}
\label{sec-2.1}

Throughout, $\nu\ll\beta$ is assumed.

In order to explain our approach,  let us consider a general
setting when ($X,Y$) is Markov process with paths from the
Skorokhod space $\mathbb{D}=\mathbb{D}_{[0,\infty)}(\mathbb{R}^2)$
of right continuous functions having limits to the left functions. Moreover,
the signal component $X$ is Markov
process itself.

\medskip
We introduce a measurable space $(\mathbb{D},\mathscr{D})$, where
$\mathscr{D}=\sigma\{(x_s,y_s),s\ge 0\}$ is the Borel
$\sigma$-algebra on $\mathbb{D}$. Let $D=(\mathscr{D}_t)_{t\ge 0}$
be the filtration of $\mathscr{D}_t=\sigma\{(x_s,y_s),s\le t\}$
and let $D^y=(\mathscr{D}^y_t)_{t\ge 0}$ be the filtration of
$\mathscr{D}^y_t=\sigma\{y_s,s\le t\}$.

As before, we write $(X^\nu_t,Y^\nu_t)$ and $(X^\beta_t,
Y^\beta_t)$, when the distribution of $X_0$ is $\nu$ or $\beta$
respectively, meaning that both pairs are defined on the same
probability space, have the same transition semigroup, but
different initial distributions.

For a bounded measurable function $f$, we introduce
$\pi^\nu_t(f):=E(f(X^\nu_t)|\mathscr{Y}^\nu_{[0,t]})$ and $
\pi^\beta_t(f):=E(f(X^\beta_t)|\mathscr{Y}^\beta_{[0,t]}). $ Since
$\pi^\nu_t(f)$ and $\pi^\beta_t(f)$ are $\mathscr{Y}^\nu_{[0,t]}$-
and $\mathscr{Y}^\beta_{[0,t]}$-measurable random variables
respectively, it is convenient to identify $\pi^\nu_t(f)$ and
$\pi^\beta_t(f)$ with some $\mathscr{D}^y_t$-measurable
functionals of trajectories $Y^\nu_{[0,t]}=\{Y^\nu_s,s\le t\}$ and
$Y^\beta_{[0,t]}=\{Y^\beta_s,s\le t\}$.

For this purpose, let $Q^\nu$ and $Q^\beta$ denote the
distributions of $(X^\nu,Y^\nu)$ and  $(X^\beta,Y^\beta)$ on
$(\mathbb{D},\mathscr{D})$ respectively and $Q^\nu_t$ and $Q^\beta_t$
be their
restrictions on $[0,t]$, so that $Q^\nu_0, Q^\beta_0$ are the
distributions of $(X^\nu_0,Y^\nu_0), (X^\beta_0,Y^\beta_0)$. We
also assume that
\begin{equation}\label{0/0}
\frac{dQ^\nu_0}{dQ^\beta_0}(x,y)=\frac{d\nu}{d\beta}(x_0).
\end{equation}
Since $(X^\nu_t,Y^\nu_t)$ and $(X^\beta_t,Y^\beta_t)$ have the
same transition law, we have $ Q^\nu\ll Q^\beta $ with
\[
\frac{dQ^\nu}{dQ^\beta}(x,y)=\frac{d\nu}{d\beta}(x_0).
\]

Without loss of generality, we assume that the filtrations $D$
and $D^y$ satisfy the general conditions with respect to
$(Q^\nu+Q^\beta)/2$.

For fixed $t$, let $H^\beta_t(y)$ be $\mathscr{D}^y_t$-measurable
functional so that $ H^\beta_t(Y^\beta)=\pi^\beta_t(f), \
\text{a.s.} $ Moreover, due to $Q^\nu\ll Q^\beta$, a version of
$H^\beta_t(y)$ can be chosen such that
the random variable
$H^\beta_t(Y^\nu)$ is well defined. Then, we identify
$\pi^{\beta\nu}_t(f)$ with $H^\beta_t(Y^\nu)$.

We do not assume that $\beta\ll\nu$ (and thus
$Q^\beta\not\ll Q^\nu$), so this construction fails for
$\pi^{\nu\beta}_t(f)$. Nevertheless, a version of $H^\nu_t(y)$ can
be chosen such that $H^\nu_t(Y^\nu)=\pi^\nu_t(f)$ a.s. and used
for the definition of $\pi^{\nu\beta}_t(f)$. Indeed, let
$\overline{Q}^\beta$ and $\overline{Q}^\nu$ be the distributions
of $Y^\nu$ and $Y^\beta$ respectively, i.e., the marginal
distributions of $Q^\beta$ and  $Q^\nu$. Obviously,
$\overline{Q}^\nu\ll\overline{Q}^\beta$ as well as
$\overline{Q}^\nu_t\ll\overline{Q}^\beta_t$; the restrictions of
$\overline{Q}^\nu$ and $\overline{Q}^\beta$ on the interval
$[0,t]$. Moreover, $
\frac{d\overline{Q}^\nu_t}{d\overline{Q}^\beta_t}(Y^\beta)=E\big(
\frac{d\nu}{d\beta}(X^\beta_0)\big|\mathscr{Y}^\beta_{[0,t]}\big).
$ Now define
\[
\pi^{\nu\beta}_t(f):=H^\nu_t(Y^\beta)
I\Bigg(\frac{d\overline{Q}^\nu_t}{d\overline{Q}^\beta_t}(Y^\beta)>0\Bigg).
\]

We introduce the decreasing filtration $
\mathscr{X}^\beta_{[t,\infty)}=\sigma\{X^\beta_s, s\ge t\}, $ the
tail $\sigma$-algebra
\begin{equation}\label{tail}
\mathscr{T}(X^\beta)=\bigcap_{t\ge 0}
\mathscr{X}^\beta_{[t,\infty)},
\end{equation}
and  $\sigma$-algebras $\mathscr{X}^\beta_t=\sigma\{X^\beta_t\},
\quad \mathscr{Y}^\beta_{[0,\infty)} = \bigvee_{t\ge 0}
\mathscr{Y}^\beta_{[0,t]}$.

Set
\begin{equation}\label{beta0}
\pi^{\beta_0}_t(f)=E\big(f(X^\beta_t)|\mathscr{Y}^\beta_{[0,t]}
\vee\mathscr{X}^\beta_0 \big).
\end{equation}

\medskip
\subsection{Filter stability}
\label{sec-2.2}
For bounded and measurable $f$, the estimate $\pi^\nu_t(f)$  is
asymptotically stable with respect to $\beta$, if
\begin{equation}
\label{q}
\lim_{t\to\infty}E\big|\pi_t^\nu(f)-\pi_t^{\beta\nu}(f)\big|=0.
\end{equation}
Note that, when the signal process takes values in a finite
alphabet and \eqref{q} holds for any bounded $f$, then \eqref{q}
and \eqref{sense'} are equivalent.

We establish below that \eqref{q} holds,
if for large values of $t$ the additional measurement $X^\beta_0$
is useless for estimation of $f(X^\beta_t)$ via $Y^\beta_{[0,t]}$
or, analogously, if the additional measurement $X^\beta_t$ is
useless for estimation of $\frac{d\nu}{d\beta}(X^\beta_0)$ via
$Y^\beta_{[0,\infty)}$.

\begin{proposition}\label{pro-1.1}
Assume $\nu\ll\beta$. Then, any of the conditions

{\rm 1.}
\begin{equation}\label{1.6}
\lim_{t\to \infty}E\big|\pi^\beta_t(f)-\pi^{\beta_0}_t(f)\big|=0,
\end{equation}

{\rm 2.}
\begin{equation}\label{1.6''}
E\Big(\frac{d\nu}{d\beta}(X^\beta_0)\big|\mathscr{Y}^\beta_{[0,\infty)}\Big)=
\lim_{t\to\infty}
E\Big(\frac{d\nu}{d\beta}(X^\beta_0)\big|\mathscr{Y}^\beta_{[0,\infty)}\vee
\mathscr{X}^\beta_{[t,\infty)}\Big),
\end{equation}
provides \eqref{q}.
\end{proposition}

\begin{proof}
Let us first show that, under $\nu\ll\beta$, for any bounded $f$
\begin{equation}
\begin{aligned}\label{impo}
&E\big|\pi^{\beta\nu}_t(f)-\pi_t^\nu(f)\big|
\\
&=E\Big|E\Big(\frac{d\nu}{d\beta}(X_0^\beta)
\big|\mathscr{Y}^\beta_{[0,t]}\Big)E\Big(f(X^\beta_t)|\mathscr{Y}^\beta_{[0,t]}
\Big)-E\Big(\frac{d\nu}{d\beta}(X_0^\beta)
f(X^\beta_t)\big|\mathscr{Y}^\beta_{[0,t]}\Big)\Big|.
\end{aligned}
\end{equation}

Write
$$
\begin{aligned}
& E\big|\pi^{\beta\nu}_t(f)-\pi_t^\nu(f)\big|
=E\frac{d\nu}{d\beta}(X_0^\beta)
\big|\pi^{\beta}_t(f)-\pi_t^{\nu\beta}(f)\big|
\\
&=EE\Big(\frac{d\nu}{d\beta}(X_0^\beta)\big|\mathscr{Y}^\beta_{[0,t]}\Big)
\big|\pi^{\beta}_t(f)-\pi_t^{\nu\beta}(f)\big|
=E\Big|E\Big(\frac{d\nu}{d\beta}(X_0^\beta)\big|\mathscr{Y}^\beta_{[0,t]}\Big)
\big(\pi^{\beta}_t(f)-\pi_t^{\nu\beta}(f)\big)\Big|
\\
&=E\Big|E\Big(\frac{d\nu}{d\beta}(X_0^\beta)\big|\mathscr{Y}^\beta_{[0,t]}\Big)
E\Big(f(X^\beta_t)\big|\mathscr{Y}^\beta_{[0,t]}\Big) -
E\Big(\frac{d\nu}{d\beta}(X_0^\beta)\pi_t^{\nu\beta}(f)\big|
\mathscr{Y}^\beta_{[0,t]}\Big)\Big|.
\end{aligned}
$$
So, it remains to show
\begin{equation}\label{2.4'}
E\Big(\frac{d\nu}{d\beta}(X_0^\beta)\pi_t^{\nu\beta}(f)\big|\mathscr{Y}^\beta_{[0,t]}
\Big)=
E\Big(\frac{d\nu}{d\beta}(X_0^\beta)f(X^\beta_t)\big|\mathscr{Y}^\beta_{[0,t]}\Big).
\end{equation}
With $\mathscr{D}^y_t$-measurable and bounded function $\Psi_t(y)$
we get
$$
\begin{aligned}
&
E\Big\{\Psi_t(Y^\beta)E\Big(\frac{d\nu}{d\beta}(X_0^\beta)\pi_t^{\nu\beta}(f)\big|
\mathscr{Y}^\beta_{[0,t]}\Big)\Big\}
=E\Big(\Psi_t(Y^\beta)\frac{d\nu}{d\beta}(X_0^\beta)\pi_t^{\nu\beta}(f)\Big)
\\
&=E\Big(\Psi_t(Y^\nu)\pi_t^\nu(f)\Big)
=E\Big(\Psi_t(Y^\nu)f(X^\nu_t)\Big)
=E\Big(\Psi_t(Y^\beta)\frac{d\nu}{d\beta}(X^\beta_0)f(X^\beta_t)\Big),
\end{aligned}
$$
and notice that \eqref{2.4'} is valid by the arbitrariness of
$\Psi_t$.

\medskip
{\it The proof of \eqref{1.6}$\Rightarrow$\eqref{q}}. Using
\eqref{impo} and
$$
E\Big(\frac{d\nu}{d\beta}(X_0^\beta)
f(X^\beta_t)\big|\mathscr{Y}^\beta_{[0,t]}\Big)=E\Big(\frac{d\nu}{d\beta}(X_0^\beta)
\pi^{\beta_0}_t(f)\big|\mathscr{Y}^\beta_{[0,t]}\Big),
$$
 we derive
$$
\begin{aligned}
&E\big|\pi^{\beta\nu}_t(f)-\pi_t^\nu(f)\big|
=E\Big|E\Big(\frac{d\nu}{d\beta}(X_0^\beta)
\big|\mathscr{Y}^\beta_{[0,t]}\Big)\pi^\beta_t(f)
-E\Big(\frac{d\nu}{d\beta}(X_0^\beta)
\pi^{\beta_0}(f)\big|\mathscr{Y}^\beta_{[0,t]}\Big)\Big|
\\
&=E\Big|E\Big(\frac{d\nu}{d\beta}(X_0^\beta)\big(\pi^\beta_t(f)-
\pi^{\beta_0}(f)\big)\big|\mathscr{Y}^\beta_{[0,t]}\Big)\Big| \le
E\frac{d\nu}{d\beta}(X_0^\beta)
\big|\pi^{\beta}_t(f)-\pi_t^{\beta_0}(f)\big|,
\end{aligned}
$$
where the Jensen inequality has been used. Let for definiteness
$|f|\le K$ with some constant $K$. Then $\pi^\beta_t(f)$,
$\pi^{\beta_0}_t(f)$ can also be chosen such that $|\pi^\beta_t(f)|$
and $|\pi^{\beta_0}_t(f)|$ are bounded by $K$. Hence, for any
$C>0$, we have
\[
E\big|\pi^{\beta\nu}_t(f)-\pi_t^\nu(f)\big|\le
CE\big|\pi^{\beta}_t(f)-\pi_t^{\beta_0}(f)\big|+2KP\Big(\frac{d\nu}{d\beta}
(X^\beta_0)>C\Big).
\]
Therefore, $
\varlimsup_{t\to\infty}E\big|\pi^{\beta\nu}_t(f)-\pi_t^\nu(f)\big|\le
2KP\big(\frac{d\nu}{d\beta}(X^\beta_0)
>C\big)
$ and by the Chebyshev inequality
$
P\big(\frac{d\nu}{d\beta}(X^\beta_0)>C\big)\le C^{-1}\to 0, \
C\to\infty.
$

\medskip
{\it The proof of \eqref{1.6''}$\Rightarrow$\eqref{q}}. By
\eqref{impo}
$$
\begin{aligned}
& E\big|\pi^{\beta\nu}_t(f)-\pi_t^\nu(f)\big|
\\
&=E\Bigg|E\Big(f(X^\beta_t)E\Big[\frac{d\nu}{d\beta}(X_0^\beta)
\big|\mathscr{Y}^\beta_{[0,t]}\Big]
\Big|\mathscr{Y}^\beta_{[0,t]}\Big)-
E\Big(f(X^\beta_t)\frac{d\nu}{d\beta}(X^\beta_0)|\mathscr{Y}^\beta_{[0,t]}\Big)\Bigg|.
\end{aligned}
$$
Notice also
$$
E\Big(f(X^\beta_t)\frac{d\nu}{d\beta}(X^\beta_0)|\mathscr{Y}^\beta_{[0,t]}\Big)=
E\Big(f(X^\beta_t)E\Big[\frac{d\nu}{d\beta}(X^\beta_0)\big|
\mathscr{Y}^\beta_{[0,\infty)}\vee
\mathscr{X}^\beta_{[t,\infty)}\Big]
\big|\mathscr{Y}^\beta_{[0,t]}\Big).
$$
Since $|f|\le K$, by the Jensen inequality we have
\begin{multline}\label{alik}
E\big|\pi^{\beta\nu}_t(f)-\pi_t^\nu(f)\big|
\\
\le
KE\Big|E\Big(\frac{d\nu}{d\beta}(X_0^\beta)\big|\mathscr{Y}^\beta_{[0,t]}\Big)
-
E\Big(\frac{d\nu}{d\beta}(X^\beta_0)\big|\mathscr{Y}^\beta_{[0,\infty)}
\vee \mathscr{X}^\beta_{[t,\infty)}\Big)\Big|.
\end{multline}
Both random processes $
E\big(\frac{d\nu}{d\beta}(X_0^\beta)\big|\mathscr{Y}^\beta_{[0,t]}\big)
$ and $
E\big(\frac{d\nu}{d\beta}(X^\beta_0)\big|\mathscr{Y}^\beta_{[0,\infty)}
\vee \mathscr{X}^\beta_{[t,\infty)}\big) $ are uniformly
integrable forward and backward martingales with respect to the
filtrations $(\mathscr{Y}^\beta_{[0,t]})_{t\ge 0}$ and
$(\mathscr{Y}^\beta_{[0,\infty)}\vee
\mathscr{X}^\beta_{[t,\infty)})_{t\ge 0}$. Therefore, they admit
limits a.s. in $t\to\infty$: $
E\big(\frac{d\nu}{d\beta}(X_0^\beta)\big|\mathscr{Y}^\beta_{[0,\infty)}\big)
$ and $ \lim\limits_{t\to\infty}
E\big(\frac{d\nu}{d\beta}(X^\beta_0)\big|\mathscr{Y}^\beta_{[0,\infty)}
\vee \mathscr{X}^\beta_{[t,\infty)}\big)$, respectively. By
\eqref{1.6''}
$$
\lim_{t\to\infty}\Big|E\Big(\frac{d\nu}{d\beta}(X_0^\beta)\big|
\mathscr{Y}^\beta_{[0,t]}\Big) -
E\Big(\frac{d\nu}{d\beta}(X^\beta_0)\big|\mathscr{Y}^\beta_{[0,\infty)}
\vee \mathscr{X}^\beta_{[t,\infty)}\Big)\Big|=0.
$$
We show also that
\begin{equation}\label{also}
\lim_{t\to\infty}E\Big|E\Big(\frac{d\nu}{d\beta}(X_0^\beta)\big|\mathscr{Y}
^\beta_{[0,t]}\Big) -
E\Big(\frac{d\nu}{d\beta}(X^\beta_0)\big|\mathscr{Y}^\beta_{[0,\infty)}
\vee \mathscr{X}^\beta_{[t,\infty)} \Big)\Big|=0.
\end{equation}
Denote by $\alpha_t$ any of $
E\Big(\frac{d\nu}{d\beta}(X_0^\beta)\big|\mathscr{Y}^\beta_{[0,t]}\Big)
$ and $
E\Big(\frac{d\nu}{d\beta}(X^\beta_0)\big|\mathscr{Y}^\beta_{[0,\infty)}
\vee \mathscr{X}^\beta_{[t,\infty)}\Big) $ and
$$
\alpha_\infty=\lim\limits_{t\to\infty}\alpha_t.
$$
It is clear that \eqref{also} holds true, if $
\lim_{t\to\infty}E|\alpha_t-\alpha_\infty|=0. $ Since
$\lim\limits_{t\to\infty}\alpha_t=\alpha_\infty$, $\alpha_t\ge 0$
and $E\alpha_t\equiv E\alpha_\infty=1$, by the Scheffe theorem we
get the desired property.

Thus the right hand side of \eqref{alik} converges to zero and the
result follows.
\end{proof}

\subsection{Connection to the gap in  \cite{K}}
\label{sec-2.3}
In \cite{K},  H. Kunita studies\footnotemark \ ergodic properties
of the filtering process $\pi_t^\nu$. He considers $\pi_t^\nu$ as a
Markov process with values in the space of probability measures
and claims (in Theorem 3.3) that there exists the
unique invariant
measure being ``limit point'' of marginal
distributions of $\pi_t^\nu$, $t\nearrow
\infty$. As was later shown in
\cite{OP}, this result is the key to
the stability analysis under \eqref{Kcond}.
\footnotetext{The notations of this paper are used here.}

Below we demonstrate that the main argument, used  in the proof of
Theorem 3.3 of \cite{K}, cannot be taken for granted. We discuss
this issue in the context of Proposition \ref{pro-1.1}. Suppose
the Markov process $X$ is  ergodic in the sense of \eqref{1.7m}
and \eqref{Kcond}. It is well known that its tail $\sigma$-algebra
$\mathscr{T}(X^\beta)$ (see \eqref{tail} for definition) is empty
almost surely. It is very tempting in this case to change the
order of intersection and supremum as follows:
\begin{equation}\label{1.11'}
\bigcap_{t\ge 0} \mathscr{Y}^\beta_{[0,\infty)}\vee
\mathscr{X}^\beta_{[t,\infty)}
=\mathscr{Y}^\beta_{[0,\infty)}\vee\mathscr{T}(X^\beta), \quad
\text{a.s.}
\end{equation}
Then, the right-hand side of \eqref{1.6''} is transformed to
$$
\begin{aligned}
& \lim_{t\to\infty}E\Big(\frac{d\nu}{d\beta}(X^\beta_0)\big|
\mathscr{Y}^\beta_{[0,\infty)}\vee
\mathscr{X}^\beta_{[t,\infty)}\Big) =
E\Big(\frac{d\nu}{d\beta}(X^\beta_0)\big|\bigcap_{t\ge 0}
\Big\{\mathscr{Y}^\beta_{[0,\infty)}\vee
\mathscr{X}^\beta_{[t,\infty)}\Big\}\Big)
\\
&\quad=
E\Big(\frac{d\nu}{d\beta}(X^\beta_0)\big|\mathscr{Y}^\beta_{[0,\infty)}\vee
\mathscr{T}(X^\beta)\Big) =
E\Big(\frac{d\nu}{d\beta}(X^\beta_0)\big|\mathscr{Y}^\beta_{[0,\infty)}\Big)
\end{aligned}
$$
and \eqref{1.6''} would be correct, regardless (!) of any other
ingredients of the problem (e.g., with $\sigma=0$ in \eqref{1.1}).

\medskip
In \cite{K}, the relation of \eqref{1.11'} type plays
the key role in verification of the uniqueness for the invariant measure corresponding
to $\pi^\nu_t,t\ge 0$.
However, the validity of  \eqref{1.11'}  is far from being
obvious. According to Williams \cite{Wil}, it "...tripped up even
Kolmogorov and Wiener" (see Y. Sinai \cite[p. 837]{Sinai} for some
details). The reader can find a discussion concerning
\eqref{1.11'} in Weizs\"acker \cite{Wei}; unfortunately, the
counterexample  there is incorrect. A proper counterexample to
\eqref{1.11'}  is given in Exercise 4.12 in Williams \cite{Wil},
which, however, seems somewhat artificial in the filtering
context. It turns out that the example, considered by  Delyon and
Zeitouni in \cite{DZ} (see \cite{Kai} by Kaijser for its earlier
discrete-time version), is nothing but another case
when \eqref{1.11'} fails.

For the reader convenience, we give below a
detailed analysis of this example.

It is important to note that the counterexamples mentioned above
do not fit exactly to the setup, considered by Kunita. They merely
indicate that \eqref{1.11'} is not evident and so the
claim of Theorem 3.3 in \cite{K} remains a conjecture.

Generally, the stability of nonlinear filters for ergodic Markov
processes remains an open problem, and some results  \cite{K-91},
\cite{St}, \cite{St-91}, \cite{BBK}, \cite{B02}, \cite{B2},
\cite{OP} based on \cite{K} have to be revised.

\section{Counterexample}
\label{sec-3}

Below we give a detailed discussion of one
counterexample to \eqref{1.11'}. Consider Markov process $X$
with values in $\mathbb{S}=\{1,2,3,4\}$, with the initial
distribution $\nu$ and the transition intensities matrix
\begin{equation}\label{2.4m}
\Lambda=\begin{pmatrix}
-1 & 1 & 0 & 0
\\
0 & -1 & 1 & 0
\\
0 & 0 & -1 & 1
\\
1 & 0 & 0 & -1
\end{pmatrix}.
\end{equation}
All states of $\Lambda$ communicate and so, $X$ is, ergodic Markov
process (see e.g., \cite{Nor}) with the unique invariant measure $
\mu =\begin{pmatrix} 1/4 & 1/4 & 1/4 & 1/4
\end{pmatrix}.
$ Let $ h(x)=I(x=1)+I(x=3), $ that is,
$$
Y_t = \int_0^t \big[I(X_s=2)+I(X_s=3)\big]ds + \sigma W_t.
$$

By Theorem \ref{theo-2.0} below, the filter is stable in this case
for any $\sigma>0$.

\subsection{Noiseless observation}
\label{sec-3.1}
Consider the case $\sigma=0$.

It will be convenient to redefine the observation
process as follows:
$$
Y_t=[I(X_t=1)+I(X_t=3)].
$$
We assume $\nu\ll\beta$ and notice that
\eqref{0/0} holds true. We omit the superscripts $\nu$ and
$\beta$, when the initial condition does not play a significant
role.
Since $X$ is ergodic Markov process, satisfying \eqref{Kcond},
$\mathscr{T}(X)= (\Omega,\varnothing)$, a.s.

\begin{proposition}\label{theo-A.1}
\begin{equation}\label{nuka}
\bigcap_{t\ge 0}\Big(\mathscr{Y}_{[0,\infty)}\vee
\mathscr{X}_{[t,\infty)}\Big) \varsupsetneq
\mathscr{Y}_{[0,\infty)}, \ \text{a.s.}
\end{equation}
\end{proposition}

\begin{proof}
It suffices to show that $X_0$ is $\bigcap_{t\ge
0}\Big(\mathscr{Y}_{[0,\infty)}\vee
\mathscr{X}_{[t,\infty)}\Big)$-measurable random variable and at
the same time $X_0\notin\mathscr{Y}_{[0,\infty)}$.

The structure of matrix $\Lambda$ admits only cyclic transitions in the following order

$$
\cdots\to \{3\}\to \{4\} \to \{1\} \to \{2\}\to \{3\}\to\cdots
$$
So, since $Y$ and $X$ jump simultaneously, $X_0$ can be  recovered exactly from
the trajectory $Y_s, s\le t$ and $X_t$ for any $t>0$, i.e.,
$X_0$ is $ \mathscr{X}_t\vee\mathscr{Y}_{[0,t]}$-measurable. Owing to
$ \mathscr{X}_t\vee
\mathscr{Y}_{[0,t]}\subset \mathscr{X}_{[t, \infty)}\vee
\mathscr{Y}_{[0,\infty)}$, $X_0$
is measurable with respect to
$$
\bigcap_{t\ge 0} \Big(\mathscr{Y}_{[0,\infty)}\vee
\mathscr{X}_{[t, \infty)}\Big).
$$

Denote by $(\tau_i)_{i\ge 1}$ the time moments  where $Y$
jumps. It is not hard
to check that $(\tau_i)_{i\ge 0}$ is independent of $(X_0,Y_0)$ and moreover
$$
\mathscr{Y}_{[0,t]}=\bigvee_{i\ge 0}\sigma\{\tau_i\le t\}\vee
\sigma\{Y_0\}.
$$

Thus for any $t\ge 0$
\begin{eqnarray}\label{valid}
\begin{aligned}
P\big(X_0=1|\mathscr{Y}_{[0,t]}\big) &= P\left(X_0=1|
\bigvee_{i\ge 0}\sigma\{\tau_i\le t\}\vee \sigma\{Y_0\}\right) \\
&=P\big(X_0=1|Y_0\big) =\frac{\nu_1}{\nu_1+\nu_3}Y_0.
\end{aligned}
\end{eqnarray}
Since \eqref{valid} is valid for any $t\ge 0$, we conclude that
$$
P\big(X_0=1|\mathscr{Y}_{[0,\infty)}\big)=\frac{\nu_1}{\nu_1+\nu_3}Y_0.
$$
Obviously $I(X_0=1)\ne \frac{\nu_1}{\nu_1+\nu_3}Y_0$ and thus
$X_0$ is not $\mathscr{Y}_{[0,\infty)}$-measurable.
\end{proof}

\subsection{Invariant measures of $\b{\pi_t}$ and the filter
instability}
\label{sec-3.2}

Since $I_t(2)+I_t(4)=1-Y_t$ and $I_t(1)+I_t(3)=Y_t$, only $I_t(1)$ and $I_t(2)$
have to be
filtered while $\pi_t(3)=Y_t-\pi_t(1)$ and
$\pi_t(4)=(1-Y_t)-\pi_t(2)$. The derivation of the filtering equations
is sketched in the appendix.

\medskip
\begin{proposition}\label{theo-A.2}
The optimal filtering estimate satisfies
$$
\begin{aligned}
d\pi_t(1)&=\big(1-\pi_{t-}(2)\big)(1-Y_{t-})dY_t+\pi_{t-}(1)Y_{t-}dY_t,
\\
d\pi_t(2)&=-\pi_{t-}(2)(1-Y_{t-})dY_t-\pi_{t-}(1)Y_{t-}dY_t
\end{aligned}
$$
subject to $\pi_0(1)=\frac{\nu_1}{\nu_1+\nu_3}Y_0$,
$\pi_0(2)=\frac{\nu_2}{\nu_2+\nu_4}(1-Y_0)$.
\end{proposition}

Let us examine the behavior of the filter from Proposition
\ref{theo-A.2}. A pair of typical trajectories are given
in Table \ref{table1} (for $Y_0=1$) and Table \ref{table2} (for
$Y_0=0$).

\begin{table}[ht]
\caption{ Typical trajectory of $\pi_t$ for $Y_0=1$.}

\begin{center} \footnotesize
\begin{tabular}{|c|c|c|c|c|c|c|}
\hline
$t$ & $[0,\tau_1)$ & $[\tau_1, \tau_2)$ & $[\tau_2,\tau_3)$ & $[\tau_3,\tau_4)$ & $[\tau_4,\tau_5)$ & ...\\
\hline
$Y_t$ & 1 & 0 & 1 & 0 & 1 & ...\\
\hline
$\pi_t(1)$ & $\frac{\nu_1}{\nu_1+\nu_3}$ & 0 & $\frac{\nu_3}{\nu_1+\nu_3}$ & 0 & $\frac{\nu_1}{\nu_1+\nu_3}$ & ... \\
\hline
$\pi_t(2)$ & 0 & $\frac{\nu_1}{\nu_1+\nu_3}$ & 0 & $\frac{\nu_3}{\nu_1+\nu_3}$ & 0 & ...\\
\hline
\end{tabular}
\end{center}
\label{table1}
\end{table}
\begin{table}[ht]
\caption{ Typical trajectory of $\pi_t$ for $Y_0=0$.}

\begin{center} \footnotesize
\begin{tabular}{|c|c|c|c|c|c|c|}
\hline
$t$ & $[0,\tau_1)$ & $[\tau_1, \tau_2)$ & $[\tau_2,\tau_3)$ & $[\tau_3,\tau_4)$
& $[\tau_4,\tau_5)$ & ...\\
\hline
$Y_t$ & 0 & 1 & 0 & 1 & 0 & ...\\
\hline
$\pi_t(1)$ & 0 & $\frac{\nu_2}{\nu_2+\nu_4}$ & 0 & $\frac{\nu_4}{\nu_2+\nu_4}$
& 0 & ... \\
\hline
$\pi_t(2)$ &  $\frac{\nu_2}{\nu_2+\nu_4}$ & 0 & $\frac{\nu_4}{\nu_2+\nu_4}$
& 0 & $\frac{\nu_2}{\nu_2+\nu_4}$& ...\\
\hline
\end{tabular}
\end{center}
\label{table2}
\end{table}

It is not hard to see that $Y$ is itself Markov chain with
values in $\{0,1\}$ and the transition intensities matrix $
\left(\begin{smallmatrix}
-1 & 1\\
1 & -1
\end{smallmatrix}
\right) $ and thus its invariant measure is $\mu'=\begin{pmatrix}
1/2 & 1/2
\end{pmatrix}$. Hence, the
invariant measure $\Phi$ of the filtering process $(\pi_t(1),\pi_t(2))$ is
concentrated on eight vectors
\begin{eqnarray*}
\begin{aligned}
\phi_1=\begin{pmatrix}
\frac{\nu_1}{\nu_1+\nu_3} \\
0
\end{pmatrix}, \quad
\phi_2=\begin{pmatrix} 0 \\
\frac{\nu_1}{\nu_1+\nu_3}
\end{pmatrix}, \quad
\phi_3=\begin{pmatrix}
\frac{\nu_3}{\nu_1+\nu_3} \\
0
\end{pmatrix}, \quad
\phi_4=\begin{pmatrix} 0 \\
\frac{\nu_3}{\nu_1+\nu_3}
\end{pmatrix}, \\
\phi_5=\begin{pmatrix}
\frac{\nu_2}{\nu_2+\nu_4} \\
0
\end{pmatrix}, \quad
\phi_6=\begin{pmatrix} 0 \\
\frac{\nu_2}{\nu_2+\nu_4}
\end{pmatrix}, \quad
\phi_7=\begin{pmatrix}
\frac{\nu_4}{\nu_2+\nu_4} \\
0
\end{pmatrix}, \quad
\phi_8=\begin{pmatrix} 0 \\
\frac{\nu_4}{\nu_2+\nu_4}
\end{pmatrix}
\end{aligned}
\end{eqnarray*}
with
\begin{eqnarray*}
\begin{aligned}
\Phi(\phi_i)=(\nu_1+\nu_3)/4, \quad i=1,2,3,4, \\
\Phi(\phi_i)=(\nu_2+\nu_4)/4, \quad i=5,6,7,8,
\end{aligned}
\end{eqnarray*}
and, consequently,  $\Phi$ is not unique. Moreover, the optimal
filter is not stable in the sense \eqref{sense'}. In fact, for
different initial conditions, the filtering distribution $\pi_t, t>0$
can ``sit" on different vectors!

\section{Main results} \label{sec-4}

\subsection{Ergodic case}
\label{sec-4.1}
Markov chain $X$ is ergodic, if and only if all entries of its
transition intensities matrix $\Lambda$  {\em communicate}, i.e.,
for any pair of indices $i$ and $j$, a string of indices
$\{\ell_1, \ldots,\ell_m\}$ can be found so that
$\lambda_{i\ell_1}\lambda_{\ell_1\ell_2}\ldots\lambda_{\ell_m j}\ne
0$ (see, e.g., \cite{Nor}). In this case, the distribution of $X_t$
converges to the positive invariant distribution $\mu$  being the
unique solution of $\Lambda^*\mu =0$ in the class of vectors with positive entries
the sum of which is equal to one.

\medskip
\begin{theorem}\label{theo-2.0}
If all states of $\Lambda$ communicate, then there exists a
positive constant $c$ such for any $\nu$ and $\beta$
$$
\varlimsup_{t\to\infty}\frac{1}{t}\log\|\pi^{\beta\nu}_t-\pi^\nu_t\|<-c,
\ \text{a.s.}
$$
\end{theorem}

\begin{remark}
Clearly, Theorem \ref{theo-2.0} provides \eqref{sense'}. Also it
allows to conclude that $
\lim\limits_{t\to\infty}\|\pi^{\beta\nu}_t-\pi^\nu_t\|=0, \
\text{a.s.} $ for $\beta$ concentrated in a single state of
\ $\mathbb{S}$. Then, in particular, we have
\[
\lim\limits_{t\to\infty}\|\pi^{\mu_0}_t-\pi^\mu_t\|=0
\]
which is the main argument in the proof of existence of the unique
invariant measure for the process $(\pi_t)_{t\ge 0}.$ This fact
corroborates Kunita's result from {\rm \cite{K}} in the finite
state space setup of Theorem \ref{theo-2.0}.
\end{remark}

\medskip
Actually, Theorem \ref{theo-2.0} verifies the logarithmic rate in
$t\to\infty$ which is in general a function of $\Lambda$, $h$ and
$\sigma$. However stronger assumptions on $\Lambda$ guarantee
exponential or logarithmic rates, regardless of $h$ and $\sigma$ ($\sigma$ is only
required to be
nonzero).

\medskip
\begin{theorem}\label{theo-2.1}
Assume  all states of $\Lambda$ communicate. Then
\begin{equation}\label{uspex2}
\varlimsup_{t\to\infty}\frac{1}{t}\log\|\pi^{\beta\nu}_t-\pi^\nu_t\|
\le -\sum_{r=1}^n\mu_r\min_{i\ne r}\lambda_{ri}.
\end{equation}
\end{theorem}
\begin{remark}
The bound \eqref{uspex2} is negative if at least one row of
$\Lambda$ has all nonzero entries.
\end{remark}

\begin{theorem}\label{theo-2.20}
Assume all entries of $\Lambda$ are nonzero.

{\rm 1.} If $\nu\ll\beta$, then
\begin{equation}\label{4.2a}
E\|\pi^{\beta\nu}_t-\pi^\nu_t\|\le
n\sum_{j=1}^n\frac{d\nu}{d\beta}(a_j) \exp\Big(-2t\min_{p\ne
q}\sqrt{\lambda_{pq}\lambda_{qp}}\Big), \ t>0.
\end{equation}

{\rm 2.} If $\nu\sim\beta$, then
\begin{equation}\label{4.2b}
\|\pi^{\beta\nu}_t-\pi^{\nu}_t\| \le
 n^2\max_j\frac{d\nu}{d\beta}(a_j)\max_j\frac{d\beta}{d\nu}(a_j)
\exp\Big(-2t\min_{p\ne q}\sqrt{\lambda_{pq}\lambda_{qp}}\Big),
 \ t>0.
\end{equation}
\end{theorem}

\subsection{Nonergodic case}
\label{sec-4.2}
Let $m\ge 2$ and $\Lambda$ be given in \eqref{LA}.
If $X_0 \in \mathbb{S}_j$, then $X$ is Markov process with values
in $\mathbb{S}_j$ with transition intensities matrix $\Lambda_j$.
We denote this process by $X^j$. In addition to $h $, introduce
column vectors $ h _j$, $j=1,\ldots ,m$ with entries $
h(a^j_1),\ldots,h(a^j_{n_j}) $ respectively.

\begin{theorem}\label{theo-2.2}
Assume the following.
\begin{enumerate}
\renewcommand{\theenumi}{A-{\rm \arabic{enumi}}}
\item  \label{A-1} For any $j$, all states of $\Lambda_j$
communicate.

\item  \label{A-2} For each $j,k$ with $j\neq k$ either
\begin{eqnarray*}
&&  h ^*_j\mu^j  \neq  h ^*_k\mu^k
\\
&&\text{or}
\\
&&  h ^*_j\diag(\mu^j)\Lambda^q_j h _j\neq
 h ^*_k\diag(\mu^k)\Lambda^q_k h _k, \ \text{for
some\ } \ 0 \le q \le n_j+n_k-1.
\end{eqnarray*}
\end{enumerate}
Then the asymptotic stability \eqref{sense'} holds true.
\end{theorem}

The condition \eqref{A-1} is inherited from Theorem \ref{theo-2.0}
to ensure the stability within each ergodic class, while under
\eqref{A-2} $\mathscr{Y}_{[0,\infty)}$ completely identifies the
class in which $X$ actually resides.

\section{Proofs for the ergodic case}
\label{sec-5}
Recall that under $m=1$, $X$ is a homogeneous ergodic Markov chain
with values in the finite alphabet $\mathbb{S}=\{a_1,\ldots,a_n\}$
with the transition intensities matrix $\Lambda$. The unique invariant
measure $\mu=(\mu_1,\ldots,\mu_n)$ is the positive
distribution on $\mathbb{S}$. Let $\nu$ be the distribution of
$X_0$ and $\beta$ a probability measure on $\mathbb{S}$. The
observation process $Y$ is defined in \eqref{1.1}. Recall that
the entries of $\pi^\nu_t$ and $\pi^{\beta\nu}_t$ are the true and
``wrong"  conditional probabilities respectively as defined in the
introduction.

\subsection{The proof of Theorem \ref{theo-2.0}}
\label{sec-5.1}
We use the method proposed by Atar and Zeitouni in \cite{AZ2}, which
is elaborated for the considered filtering setup  for reader
convenience.

Recall the following facts from the theory of nonnegative
matrices. For a pair $(p,q)$ of nonnegative measures on
$\mathbb{S}$ (i.e., vectors with nonnegative entries), the Hilbert
projective metric $H(p,q)$ is defined as the following (see, e.g., \cite{Seneta}):
\begin{equation}\label{hilbert}
H(p,q)= \left\{
\begin{array}{ll}
\log \frac{\max\limits_{j:q_j>0}(p_j/q_j)}
{\min\limits_{i:q_i>0}(p_i/q_i)}, &
 p\sim q, \\
\infty, & p\not \sim q.
\end{array}
 \right.
\end{equation}
The Hilbert metric is known to satisfy the following properties:
\begin{enumerate}
\item \label{prop1}
$H(c_1p,c_2q)=H(p,q)$ for any positive constants $c_1$ and $c_2$.

\item \label{prop3} for matrix $A$ with nonnegative entries
($A_{ij})$,
$$
H\big(Ap,Aq\big)\le \tau(A)H\big(p,q\big) \quad\text{(see, e.g.,
\cite{Seneta})}
$$
where $ \tau(A) =\frac{1-\sqrt{\psi(A)}}{1+\sqrt{\psi(A)}} $ is
the Birkhoff contraction coefficient with
$$ \psi(A)=
\min\limits_{i,j,k,\ell}\frac{A_{ik}A_{j\ell}}{A_{i\ell}A_{jk}}.
$$

\item \label{prop2} $ \|p-q\|\le \frac{2}{\log 3}H(p,q)
\quad\text{(\cite[Lemma 1]{AZ2})}. $
\end{enumerate}

\medskip
Returning to the filtering problem, let us first consider the
special case when $\nu=\mu$ and thus the signal $X^\mu$ is the stationary Markov
chain. It is well known that $
\pi^\mu_t=\eta^\mu_t/\langle \mathbf{1}, \eta^\mu_t\rangle, $
where $\mathbf{1}$ denotes the vector with unit entries, $\langle
\cdot, \cdot \rangle$ is the usual inner product and $\eta^\mu_t$
solves the Zakai equation
\begin{equation}\label{zakai}
d\eta^\mu_t = \Lambda^* \eta^\mu_t dt + \sigma^{-2} \diag(h)
\eta^\mu_t dY^\mu_t
\end{equation}
subject to $\eta^\mu_0 = \mu$. Similarly, $
\pi^{\beta\mu}_t=\eta^{\beta\mu}_t/\langle \mathbf{1},
\eta^{\beta\mu}_t\rangle $, where $\eta^{\beta\mu}_t$ is the
solution of \eqref{zakai} subject to $\eta^{\beta\mu}_0=\beta$.

The Zakai equation  possesses the unique strong solution which is
linear with respect to the initial condition. Hence,
$\eta^\mu_t=J_{[0,t]}\mu$ and $\eta^{\beta\mu}_t=J_{[0,t]}\beta$,
$t>0$, where $J_{[0,t]}$ is the random Cauchy matrix corresponding
to \eqref{zakai}.

The  matrix $J_{[0,t]}$ can be factored (here $\lfloor t\rfloor$
is the integer part of $t$):
$$
J_{[0,t]} = J_{[\lfloor t\rfloor,t]}\left(\prod_{n=2}^{\lfloor
t\rfloor}J_{[n-1,n]}\right)J_{[0,1]}.
$$
The properties of Hilbert metric, listed above,  provide
\begin{align*}
& \big\|\pi^{\mu}_t - \pi^{\beta\mu}_t\big\| \le \frac{2}{\log 3}
H\big(\pi^\mu_t, \pi^{\beta\mu}_t\big) = \frac{2}{\log
3}H\big(J_{[0,t]}\mu, J_{[0,t]}\beta\big)
\\
&\le \frac{2}{\log 3}\tau\big(J_{[\lfloor t\rfloor,t]}\big)
\prod_{n=2}^{\lfloor t\rfloor}
\tau\big(J_{[n-1,n]}\big)H\big(J_{[0,1]}\mu, J_{[0,1]}\beta\big).
\end{align*}

Assume for a moment that $H\big(J_{[0,1]}\mu, J_{[0,1]}\beta\big)<\infty$ a.s.
Then
\begin{multline}
\label{limlim} \varlimsup_{t\to\infty}\frac{1}{t}\log
\big\|\pi^{\mu}_t - \pi^{\beta\mu}_t\big\| \le
\varlimsup_{t\to\infty}\frac{1}{\lfloor t \rfloor}
\sum_{n=2}^{\lfloor t \rfloor} \log \tau\big(J_{[n-1,n]}\big)
\\
\le \varlimsup_{t\to\infty}\frac{1}{\lfloor t\rfloor}
\sum_{n=2}^{\lfloor t \rfloor}\big\{-1\vee \log
\tau\big(J_{[n-1,n]}\big)\big\}=
E\big[-1\vee\log\tau\big(J_{[0,1]})\big]\le 0.
\end{multline}

The equality is implied by the law of large numbers, which
is valid since $-1\le \big\{-1\vee\log
\tau\big(J_{[n-1,n]}\big)\big\} \le 0$ and
$\log\tau\big(J_{[n-1,n]}\big)$ is generated by
$$
\{X^\mu_s-X^\mu_{n-1}, \ W_s-W_{n-1}\}, \ n-1\le s<n,
$$
where the processes $X^\mu$ and $W$ are independent and $X^\mu$ is an
ergodic Markov chain.

Let $J^\nu_{[n-1,n]}$ be the matrices defined similarly to
$J_{[n-1,n]}$ with $Y^\mu$ replaced by $Y^\nu$. Recall that $\mu$
is the positive measure on $\mathbb{S}$, so that $\nu\ll\mu$ and,
in turn, $\overline{Q}^\nu\ll \overline{Q}^\mu$ (here
$\overline{Q}^\mu$ is the distribution of $Y^\mu$).

Since \eqref{limlim} holds $\overline{Q}^\mu$-a.s., it also holds
$\overline{Q}^\nu$-a.s.,  i.e., with $J_{[n-1,n]}$ replaced by
$J^\nu_{[n-1,n]}$ which gives

\begin{theorem}\label{theo-az} {\rm (version of Theorem 1(a) in Atar and Zeitouni,
\cite{AZ2})} Assume that all states of $\Lambda$ communicate, i.e.,
$X$ is an ergodic Markov chain. Assume $J_{[0,1]}\beta$ and
$J_{[0,1]}\nu$ have positive entries a.s. Then,
\begin{equation}
\label{azthm} \varlimsup_{t\to\infty}\frac{1}{t}\log
\big\|\pi^{\nu}_t - \pi^{\beta\nu}_t\big\| \le E \big[-1\vee\log
\tau\big(J_{[0,1]}\big)\big], \ \text{a.s.}.
\end{equation}
\end{theorem}

Now the statement of Theorem \ref{theo-2.0} follows from the lemma below.
\begin{lemma}\label{lem-5.2m}
The right-hand side of \eqref{azthm} is strictly negative.
\end{lemma}

\begin{proof}
It suffices to show that all entries of $J_{[0,1]}$ are positive
a.s. For fixed $i,j$, we have
$$
J_{[0,t]}(i,j)=\delta_{ij}+\int_0^tJ_{[0,s]}(i,j)\big[\lambda_{ii}ds+
\sigma^{-2}h(a_i)dY^\mu_s\big] +\int_0^t\sum_{r\ne
i}\lambda_{ri}J_{[0,s]}(r,j)ds.
$$

With the help of It\^o formula and with
$$
\phi_t(i)=\exp\big\{\lambda_{ii}t+\sigma^{-2}h(a_i)Y^\mu_t-(1/2)\sigma^{-2}h^2(a_i)
t\big\}
$$
we derive
\begin{equation}\label{5.2x}
\begin{aligned}
J_{[0,t]}(j,j)&=\phi_t(j)\Big(1+\int_0^t\phi^{-1}_s(j)\sum_{r\ne
j}\lambda_{rj} J_{[0,s]}(r,j)ds\Big),
\\
J_{[0,t]}(i,j)&=\phi_t(i)\int_0^t\phi^{-1}_s(i)\sum_{r\ne
i}\lambda_{ri} J_{[0,s]}(r,j)ds, \quad i\ne j.
\end{aligned}
\end{equation}
Also notice that the entries of $J_{[0,t]}$ are unnormalized
conditional probabilities and so nonnegative a.s. Since all
states of $\Lambda$ communicate, for pair of indices $(i,j)$ there
is a string of indexes $j=i_\ell,\ldots,i_1=i$ such that
$\lambda_{i_\ell i_{\ell-1}}\ldots\lambda_{i_2 i_1}>0$. So from
\eqref{5.2x}, it follows that a.s.
\begin{align*}
&
J_{[0,t]}(i_\ell,i_\ell)\ge \phi_t(i_\ell)>0,
\\
& J_{[0,t]}(i_{\ell-1},i_{\ell})\ge
\phi_t(i_{\ell-1})\int_0^t\phi^{-1}_s(i_{\ell-1})\lambda_{i_\ell
i_{\ell-1}} J_{[0,s]}(i_\ell,i_\ell)ds>0,
\\
&J_{[0,t]}(i_{\ell-2},i_{\ell})\ge
\phi_t(i_{\ell-2})\int_0^t\phi^{-1}_s(i_{\ell-2})\lambda_{i_{\ell-1}i_{\ell-2}}
J_{[0,s]}(i_{\ell-1},i_\ell)ds>0
\end{align*}
for any $t>0$, and so on until we get $J_{[0,t]}(i_1,i_\ell)>0$,
$t>0$.
\end{proof}

\subsection{The proof of Theorem \ref{theo-2.1}}
\label{sec-5.2}
Denote
$
\rho_{ji}(t)=P\big(X^\beta_0=a_j|\mathscr{Y}^\beta_{[0,t]},
X^\beta_t=a_i\big).
$
If $\beta$ is a positive distribution, then by
Lemma 9.5 in \cite[Chap. 9]{LSI} we have
\begin{equation}\label{smootheq}
\begin{split}
&\rho_{ji}(0)=
  \begin{cases}
    1, & j=i \\
    0, & j\ne i
  \end{cases}
\\
&\frac{d\rho_{ji}(t)}{dt} = \sum_{r\neq i}\frac{\lambda_{ri}
\pi^\beta_t(r)}{\pi^\beta_t(i)}\big(\rho_{jr}(t)-\rho_{ji}(t)\big),
\quad i=1,\ldots,n.
\end{split}
\end{equation}

\begin{remark}\label{rem-3}
By the arguments, used in the proof of Lemma \ref{lem-5.2m}, it
can be readily shown that $\pi^\beta_t(i)>0$ a.s., $i=1,\dots,n$
for any $t>0$. Then  \eqref{smootheq} remains valid
for $t>t_0$ for any $t_0>0$ initialized by
$$
\rho_{ji}(t_0)=P\big(X^\beta_0=a_j|\mathscr{Y}^\beta_{[0,t_0]},
X^\beta_{t_0}=a_i\big).
$$
\end{remark}

Set $ i^\diamond(t) =
\argmax_{i\in\mathbb{S}}\rho_{ji}(t) $ and $ i_\diamond(t)
= \argmin_{i\in\mathbb{S}}\rho_{ji}(t) $ (if the maximum or the minimum
are attained at
several indices, the lowest one is taken by convention).
 Set
\begin{equation}\label{diamond2}
\rho^\diamond(t):=\rho_{ji^\diamond(t)}(t)\quad\text{and}\quad
\rho_\diamond(t):=\rho_{ji_\diamond(t)}(t).
\end{equation}

\begin{lemma}\label{lem-!!!}
The processes $\rho^\diamond(t)$ and $\rho_\diamond(t)$ have absolutely
continuous paths with
\begin{equation}\label{5.10m}
\begin{aligned}
d\rho^\diamond(t)=\sum_{i=1}^nI(i^\diamond(t)=i)\dot{\rho}_{ji}(t)dt,
\\
d\rho_\diamond(t)=\sum_{i=1}^nI(i_\diamond(t)=i)\dot{\rho}_{ji}(t)dt.
\end{aligned}
\end{equation}
\end{lemma}

The proof of this lemma uses two results formulated in
Propositions \ref{pro-5.-4} and \ref{pro-5.4} below.

\begin{proposition}\label{pro-5.-4}{\rm (Theorem A.6.3 in
Dupuis and Ellis \cite{DupEl}).}
Let $g=g(t)$ be an absolutely continuous function mapping of
$[0,1]$ into $\mathbb{R}$. Then for each real number $a$ the set
$\{t:g(t)=a,\dot{g}(t)\neq 0\}$ has Lebesgue measure $0$.
\end{proposition}

\begin{proposition}\label{pro-5.4}
Let $X(t,\omega)$ be a random process with absolutely continuous
paths with respect to $dt$ in the
 sense that there exists a measurable random process $x(t,\omega)$ such that
$\int_0^t|x(s,\omega)|ds<\infty$ a.s., $t>0$, and
\begin{equation}\label{X(s)}
X(t,\omega)=X(0,\omega)+\int_0^tx(s,\omega)ds.
\end{equation}

Then
\[
|X(t,\omega)|=|X(0,\omega)|+\int_0^t\sign(X(s,\omega))x(s,\omega)ds,
\]
where $sign(0)=0$.
\end{proposition}
\begin{proof}
Set
$V_t(\omega)=\int_0^t|x(s,\omega)|ds$ and notice that
for any $t'\le t''$ it holds that
$$
\big||X(t'',\omega)|-|X(t',\omega)|\big|\le |X(t'',\omega)-X(t',\omega)|
\le (V_{t''}(\omega)-V_{t'}(\omega)).
$$
Hence, for fixed $\omega$, the function $|X(t,\omega)|$ possesses
bounded total variation for any finite time interval. Denote
$U_t(\omega)$ this total variation corresponding to $[0,t]$.
Obviously, $dU_t(\omega)\ll dV_t(\omega)\ll dt$. Recall that
$U_t(\omega)=U'_t(\omega)+U''_t(\omega)$, where $U'_t(\omega)$,
$U''_t(\omega)$ are increasing continuous in $t$ functions such
that for any $t>0$ and measurable set $A$ from $\mathbb{R}_+$, $
\int_{A\cap[0,t]}dU''_s(\omega)=0 $ and $
\int_{(\mathbb{R}_+\setminus A)\cap[0,t]}dU'_s(\omega)=0$,
and at the same time, $|X(t,\omega)|=U''_t(\omega)-U'_t(\omega)$. Since
$dU'_t\ll dU_t(\omega)$, $dU''_t\ll dU_t(\omega)$,
it follows
$d|X(t,\omega)|\ll dU_t(\omega)\ll dV_t(\omega)\ll dt$ and so that
\begin{equation}\label{|X|}
|X(t,\omega)|=|X(0,\omega)|+\int_0^tg(s,\omega)ds
\end{equation}
though we may not claim that $g(t,\omega)$ is measurable in $(t,\omega)$.

Now, we show that $\sign(X(s,\omega))x(s,\omega)$ is a measurable version of
$g(s,\omega)$. By \eqref{X(s)}, we have
$
X^2(t,\omega)=X^2(0,\omega)+2\int_0^tX(s,\omega)x(s,\omega)ds.
$
At the same time by \eqref{|X|} it holds
$
|X(t,\omega)|^2=|X(0,\omega)|^2+2\int_0^t|X(s,\omega)|g(s,\omega)ds.
$
Hence, the following identity is valid:
for any $t\ge 0$
$$
\int_0^t|X(s,\omega)|g(s,\omega)ds\equiv\int_0^tX(s,\omega)x(s,\omega)ds.
$$
Therefore,
$
|X(s,\omega)|g(s,\omega)=X(s,\omega)x(s,\omega)
$
for almost all $s$ with respect to Lebesgue measure.
Consequently, we have $I(|X(s,\omega)|\ne 0)g(s,\omega)=\sign(X(s,\omega))x(s,\omega)$
for almost all $s$ with respect to Lebesgue measure. It remains to show that
$$
I(X(s,\omega)=0)g(s,\omega)=0
$$
for almost all $s$ with respect to Lebesgue measure. Taking into account \eqref{|X|},
it suffices to prove that $\int_0^\infty I(X(s,\omega)=0)d|X(s,\omega)|=0$, a.s.
On the other hand, whereas $d|X(t,\omega)|\ll dV_t(\omega)$, it suffices to show that
$\int_0^\infty I(X(s,\omega)=0)dV_s(\omega)=0$, a.s.
The latter holds by Proposition \ref{pro-5.-4}.
\end{proof}

\medskip
Now we give the proof for Lemma \ref{lem-!!!}.
\begin{proof}
Let us introduce $\rho^{\diamond,i}(t)=\rho_{j1}\vee\rho_{j2}\vee\cdots\vee\rho_{ji}$
and
$\rho_{\diamond,i}(t)=\rho_{j1}\wedge\rho_{j2}\wedge\cdots\wedge\rho_{ji}$ and
notice that $\rho^{\diamond,n}(t)=\rho^\diamond(t)$,
$\rho_{\diamond,n}(t)=\rho_\diamond(t)$.

The use of obvious identities
$$
\begin{aligned}
\rho^{\diamond,2}(t)+\rho_{\diamond,2}(t)&=\rho_{j1}(t)+\rho_{j2}(t),
\\
\rho^{\diamond,2}(t)-\rho_{\diamond,2}(t)&=|\rho_{j1}(t)-\rho_{j2}(t)|
\end{aligned}
$$
and the fact, provided by Proposition \ref{pro-5.4}, that
$
d|\rho_{j1}(t)-\rho_{j2}(t)|=p(t,\omega)dt
$
with measurable derivative $p(\omega,t)$, allow us to claim
that $\rho^{\diamond,2}(t)$ and $\rho_{\diamond,2}(t)$ are absolutely continuous
with respect to $dt$ with measurable derivatives.

Further, taking into account
$\rho^{\diamond,i}(t)=\rho^{\diamond,i-1}(t)\vee\rho_{ji}$ and
$\rho_{\diamond,i}(t)=\rho_{\diamond,i-1}(t)\wedge\rho_{ji}(t)$
and consequent identities
$$
\begin{aligned}
\rho^{\diamond,i}(t)+\rho^{\diamond,i-1}(t)\wedge\rho_{ji}(t)&=\rho^{\diamond,i-1}(t)+
\rho_{ji}(t)
\\
\rho^{\diamond,i}(t)-\rho^{\diamond,i-1}(t)\wedge\rho_{ji}(t)&=|\rho^{\diamond,i-1}(t)-
\rho_{ji}(t)|
\\
\rho_{\diamond,i-1}(t)\vee\rho_{ji}(t)+\rho_{\diamond,i}(t)&=\rho_{\diamond,i-1}(t)+
\rho_{ji}(t)
\\
\rho_{\diamond,i-1}(t)\vee\rho_{ji}(t)-\rho_{\diamond,i}(t)&=|\rho_{\diamond,i-1}(t)-
\rho_{ji}(t)|
\end{aligned}
$$
absolute continuity for $\rho^\diamond(t)$ and
$\rho_\diamond(t)$ is verified by the induction method.

Thus, $d\rho^\diamond(t)=u(t)dt$ with some density $u(t)$
such that $\int_0^t|u(s)|ds<\infty$ a.s., $t>0$. On the other hand,
since $\sum_{i=1}^nI(i^\diamond(t)=i)=1$, we have
$$
\rho^\diamond(t)=\rho^\diamond(0)+\int_0^t\sum_{i=1}^nI(i^\diamond(s)=i)u(s)ds.
$$
So, it suffices to show that for any $t>0$ and any $i=1,2\ldots,n$
$$
\int_0^tI(i^\diamond(s)=i)|u(s)-\dot{\rho}_{ji}(s)|ds=0, \ \text{a.s.}
$$
The latter holds true by Proposition \ref{pro-5.-4}, since
\begin{align*}
&\int_0^tI(i^\diamond(s)=i)|u(s)-\dot{\rho}_{ji}(s)|ds \\
&=\int_0^t I(\rho^\diamond(s) -\rho_{ji}(s)=0)|u(s)-\dot{\rho}_{ji}(s)|ds  \\
&=\int_0^t I\big(\rho^\diamond(s) -\rho_{ji}(s)=0, u(s)-\dot{\rho}_{ji}(s)
\ne 0\big)|u(s)-\dot{\rho}_{ji}(s)|ds =0.
\end{align*}
\end{proof}

\medskip
\begin{lemma}\label{lem-3.1}
Under the assumptions of Theorem {\rm \ref{theo-2.1}},
\begin{eqnarray}\label{3.3b}
\varlimsup_{t\to\infty}\frac{1}{t}\log\max_{1\le j,k,\ell\le n}
\big|\rho_{jk}(t)-\rho_{j\ell}(t)\big|\le -\sum_{r=1}^n
\mu_r\min_{i\ne r}\lambda_{ri}.
\end{eqnarray}
\end{lemma}

\medskip
\begin{proof}
By \eqref{smootheq} and \eqref{5.10m}, we have
\footnote{In \eqref{2.7},\ldots,\eqref{bnd}  we use for brevity a form of differential
equalities (inequalities) which are valid for
any $\omega$ and almost all $t$ with respect to Lebesgue measure.}
\begin{equation}\label{2.7}
\begin{split}
& \frac{d\rho_\diamond(t)}{dt}=\sum_{r\neq i_\diamond(t)}
\frac{\lambda_{ri_\diamond(t) }\pi^\beta_t(r)}
{\pi^\beta_t(i_{\diamond}(t))}\big(\rho_{jr}(t)-\rho_\diamond(t)\big)
\\
& \frac{d\rho^\diamond(t)}{dt}=\sum_{r\neq i^\diamond(t)}
\frac{\lambda_{ri^\diamond(t)}\pi^\beta_t(r)}
{\pi^\beta_t(i^\diamond(t))}\big(\rho_{jr}(t)-\rho^\diamond(t)\big).
\end{split}
\end{equation}

In what follows, we will omit the time variable in  $i_\diamond(t)$
and  $i^\diamond(t)$ for
 brevity.

Set $\triangle_t=\rho^\diamond(t)-\rho_\diamond(t)$. By
\eqref{2.7} we have
\begin{eqnarray}\label{3.4a}
\frac{d\triangle_t}{dt}&=& -\sum_{r\ne
i^\diamond}\frac{\lambda_{ri^\diamond}
\pi^\beta_t(r)}{\pi^\beta_t(i^\diamond)}\big(\rho^\diamond(t)-\rho_{jr}(t)\big)
- \sum_{r\ne i_\diamond}\frac{\lambda_{ri_\diamond}
\pi^\beta_t(r)}{\pi^\beta_t(i_\diamond)}\big(\rho_{jr}(t)-\rho_\diamond(t)\big)
\nonumber\\
&=&-\triangle_t\Bigg(\frac{\lambda_{i_\diamond i^\diamond}
\pi^\beta_t(i_\diamond)}{\pi^\beta_t(i^\diamond)} +
\frac{\lambda_{i^\diamond i_\diamond}
\pi^\beta_t(i^\diamond)}{\pi^\beta_t(i_\diamond)}\Bigg)
\\
&-&\triangle_t\left(\sum_{\substack{r\ne i^\diamond(t)\\r\ne
i_\diamond(t)}} \Big[\frac{\lambda_{ri^\diamond}
\pi^\beta_t(r)}{\pi^\beta_t(i^\diamond)}\Big(\frac{\rho^\diamond(t)-\rho_{jr}(t)}
{\triangle_t}\Big) + \frac{\lambda_{ri_\diamond}
\pi^\beta_t(r)}{\pi^\beta_t(i_\diamond)}\Big(\frac{\rho_{jr}(t)-\rho_\diamond(t)}
{\triangle_t}\Big)\Big]\right). \nonumber
\end{eqnarray}
Letting $0/0=1/2$, set
$\alpha_r(t)=\frac{\rho^\diamond(t)-\rho_{jr}(t)}{\triangle_t}$.
Then, we get $ 1-\alpha_r(t)=\frac{\rho_{jr}(t)-\rho_\diamond(t)}
{\triangle_t} $ and  $0\le\alpha_r(t)\le 1$ and  \eqref{3.4a}
implies
\begin{eqnarray}\label{bnd}
\begin{aligned}
\frac{d\triangle_t}{dt}=&-\triangle_t\Bigg(\frac{\lambda_{i_\diamond
i^\diamond} \pi^\beta_t(i_\diamond)}{\pi^\beta_t(i^\diamond)} +
\frac{\lambda_{i^\diamond i_\diamond}
\pi^\beta_t(i^\diamond)}{\pi^\beta_t(i_\diamond)}\Bigg)
\\
&-\triangle_t\left(\sum_{\substack{r\ne i^\diamond(t)\\r\ne
i_\diamond(t)}} \Big[\alpha_r(t)\frac{\lambda_{ri^\diamond}
\pi^\beta_t(r)}{\pi^\beta_t(i^\diamond)} +\big(1-\alpha_r(t)\big)
\frac{\lambda_{ri_\diamond}
\pi^\beta_t(r)}{\pi^\beta_t(i_\diamond)}\Big]\right)
\\
\le& -\triangle_t\Big(\lambda_{i_\diamond i^\diamond}
\pi^\beta_t(i_\diamond) + \lambda_{i^\diamond i_\diamond}
\pi^\beta_t(i^\diamond)\Big)
\\
&-\triangle_t\left(\sum_{\substack{r\ne i^\diamond(t)\\r\ne
i_\diamond(t)}} \Big[\alpha_r(t)\lambda_{ri^\diamond}
+\big(1-\alpha_r(t)\big)\lambda_{ri_\diamond}\Big]\pi^\beta_t(r)\right)
\\
\le& -\triangle_t\left(\lambda_{i_\diamond i^\diamond}
\pi^\beta_t(i_\diamond) + \lambda_{i^\diamond i_\diamond}
\pi^\beta_t(i^\diamond) +\sum_{\substack{r\ne i^\diamond(t)\\r\ne
i_\diamond(t)}}
\Big[\lambda_{ri^\diamond}\wedge\lambda_{ri_\diamond}\Big]\pi^\beta_t(r)\right).
\end{aligned}
\end{eqnarray}
Recall that all offdiagonal entries of $\Lambda$ are nonnegative
and $\sum_{r=1}^n \lambda_{ir}=0$ for any $i$. Then,
$\big|\lambda_{i_\diamond i^\diamond}|\wedge| \lambda_{i_\diamond
i_\diamond}\big| \ge \lambda_{i_\diamond i^\diamond}, $ $
\big|\lambda_{i^\diamond i^\diamond}|\wedge|\lambda_{i^\diamond
i_\diamond}\big| \ge \lambda_{i^\diamond i_\diamond}, $ and
\eqref{bnd} provides
\begin{eqnarray*}
\begin{aligned}
\frac{d\triangle_t}{dt}\le& -\triangle_t\sum_{r=1}^n
\Big(|\lambda_{ri^\diamond}|\wedge|\lambda_{ri_\diamond}|\Big)\pi^\beta_t(r)
\le-\triangle_t\sum_{r=1}^n\min_{1\le i\le
n}|\lambda_{ri}|\pi^\beta_t(r)
\\
&= -\triangle_t\sum_{r=1}^n \pi^\beta_t(r)\min_{i\ne
r}\lambda_{ri}.
\end{aligned}
\end{eqnarray*}
Since the derivative $\frac{d\triangle_t}{dt}$ is defined for each
$\omega$ and almost everywhere (a.e.) in $t$ with respect to $dt$,
the above inequality $ \frac{d\triangle_t}{dt} \le
-\triangle_t\sum_{r=1}^n \pi^\beta_t(r)\min_{i\ne r}\lambda_{ri} $
is also valid a.e. So, it allows us to define a.e. the function
$$
H(t)=-\triangle_t\sum_{r=1}^n \pi^\beta_t(r)\min_{i\ne
r}\lambda_{ri}-\frac{d\triangle _t}{dt}.
$$
Moreover, for the definiteness, we may  redefine $H(t)$ everywhere so as
$H(t)\ge 0$. Then we have
$$
d\triangle_t=-\Big[\triangle_t\sum_{r=1}^n \pi^\beta_t(r)\min_{i\ne r}\lambda_{ri}
+H(t)\Big]dt.
$$
Notice also that $\int_0^t|H(s)|ds<\infty$, a.s. for any $t>0$ and recall that
$\triangle_0=1$. Then, we get
$$
\triangle_t=\exp\left(-\int_0^t\sum_{r=1}^n \pi^\beta_s(r)\min_{i\ne
r}\lambda_{ri}ds\right)-\int_0^t\exp\left(-\int_s^t\sum_{r=1}^n \pi^\beta_u(r)\min_{i\ne
r}\lambda_{ri}du\right)H(s)ds
$$
and in turn
$$
\frac{1}{t}\log\triangle_t\le-\sum_{r=1}^n\Big(\min\limits_{i\ne
r}\lambda_{ri}\Big) \frac{1}{t} \int_0^t\pi^\beta_s(r)ds.
$$
So, it
is left to verify that
\begin{eqnarray}
\label{leftto}
\lim_{t\to\infty}\frac{1}{t}\int_0^t\pi^\beta_s(r)ds=\mu_r,  \quad
\text{a.s.}
\end{eqnarray}
Similarly to \eqref{1.3n}, $\pi^\beta_t$ satisfies
\begin{equation*}
\begin{split}
\pi^\beta_0=&\beta
\\
d\pi^\beta_t=&\Lambda^*\pi^\beta_t dt +
\sigma^{-2}\big(\diag(\pi_t^\beta)-\pi_t^\beta
(\pi_t^\beta)^*\big) h (dY^\beta_t- h ^*\pi^\beta_tdt).
\end{split}
\end{equation*}
Recall that $\sigma^{-1}\big(Y^\beta_t-\int_0^t h
^*\pi^\beta_sds\big)$ is the innovation Wiener process (see, e.g.,
Theorem 9.1 in Chapter 10 in \cite{LSII}). Hence $
M_t=\int_0^t\big(\diag(\pi_s^\beta)-\pi_s^\beta
(\pi_s^\beta)^*\big) h (dY^\beta_s- h ^*\pi^\beta_sds) $ is
vector-valued continuous martingale. Its entries $M_t(i)$,
$i=1,\ldots,n$, have predictable quadratic variation processes
$\langle M(i)\rangle_t$ with the following property: for some positive
constant $c$, $ d\langle M(i)\rangle_t\le cdt. $ Then by Theorem
10 in Chapter 3 in \cite{LSMar}, $
\lim_{t\to\infty}\frac{1}{t}M_t(i)=0$, a.s. This fact and
the boundedness of $\pi^\beta_t$ provide $
\Lambda^*\lim_{t\to\infty}\frac{1}{t}\int_0^t\pi^\beta_sds=0. $
The vector $Z_t=\frac{1}{t}\int_0^t\pi^\beta_sds$ has nonnegative
entries, whose sum equals 1. Therefore the limit vector
$Z_\infty$, obeying the same property, is the unique solution of
the linear algebraic equation $\Lambda^*Z_\infty=0$, i.e.,
$Z_\infty=\mu$.
\end{proof}

\bigskip
To prove Theorem \ref{theo-2.1}, without loss generality, due to
Remark \ref{rem-3},
we may assume that $\nu\sim\beta$. Then, we
show that for any $t\ge 0$ and $i=1,\ldots,n$
\begin{equation}\label{5.21m}
\big|\pi_t^\nu(i)-\pi_t^{\beta\nu}(i)\big| \le
n\max_j\frac{d\nu}{d\beta}(a_j)\max_j\frac{d\beta}{d\nu}(a_j)\max_{1\le
i,j,k\le d} \big|\rho_{ji}(t)-\rho_{jk}(t)\big|.
\end{equation}

Recall that $Q^\nu$ and $Q^\beta$ are distributions of $(X^\nu,
Y^\nu)$ and $(X^\beta, Y^\beta)$ respectively, which are
equivalent, by virtue of $\nu\sim\beta$, with
$$
\frac{dQ^\beta}{dQ^\nu}(X^\nu,Y^\nu)\equiv\frac{d\beta}{d
\nu}(X^\nu_0)\quad\mbox{and} \quad
\frac{dQ^\nu}{dQ^\beta}(X^\beta,Y^\beta)\equiv\frac{d\nu}{d\beta}(X^\beta_0).
$$
Now, we show that for any $i=1,\ldots,d$ and $t>0$, $Q^\nu$- and
$Q^\beta$-a.s.
\begin{eqnarray}
\label{nuform} \pi^{\beta\nu}_t(i)=\frac
{\sum_{j=1}^n\Big(\frac{d\beta}{d\nu}(a_j)P\big(X^\nu_0=a_j,
X_t^\nu=a_i|\mathscr{Y}^\nu_{[0,t]}\Big)}{E\Big(\frac{d\beta}{d\nu}(X^\nu_0)|\mathscr{Y}^\nu_{[0,t]}\Big)}.
\end{eqnarray}
To this end, with any bounded $\mathscr{D}^y_t$-measurable
function $\psi_t(y)$, write
\begin{eqnarray*}
\begin{aligned}
& E
\psi_t(Y^\nu)\pi^{\beta\nu}_t(i)E\Big(\frac{d\beta}{d\nu}(X_0^\nu)|
\mathscr{Y}^\nu_{[0,t]}\Big) = E \psi_t(Y^\nu)
\pi^{\beta\nu}_t(i)\frac{d\beta}{d\nu}(X_0^\nu)
\\
&= E \psi_t(Y^\nu)
\pi^{\beta\nu}_t(i)\frac{dQ^\beta}{dQ^\nu}(X^\nu,Y^\nu) =
E\psi_t(Y^\beta) \pi^{\beta}_t(i)
\\
&= E\psi_t(Y^\beta) I(X^\beta_t=a_i) = E
\psi_t(Y^\nu)I(X^\nu_t=a_i)\frac{dQ^\beta}{dQ^\nu}(X^\nu,Y^\nu)
\\
&= E \psi_t(Y^\nu)I(X^\nu_t=a_i)\frac{d\beta}{d\nu}(X^\nu_0) =
E\psi_t(Y^\nu)E\Big(I(X^\nu_t=a_i)\frac{d\beta}{d\nu}(X^\nu_0)\big|\mathscr{Y}^\nu_{[0,t]}\Big).
\end{aligned}
\end{eqnarray*}
Hence, by the arbitrariness of $\psi_t(y)$,
$$
\pi^{\beta\nu}_t(i)E\Big(\frac{d\beta}{d\nu}(X_0^\nu)|\mathscr{Y}^\nu_{[0,t]}\Big)
=E\big(I(X^\nu_t=a_i)\frac{d\beta}{d\nu}(X^\nu_0)|
\mathscr{Y}^\nu_{[0,t]}\big).
$$
Further, $Q^\nu\sim Q^\beta$ provides $
E\big(\frac{d\beta}{d\nu}(X_0^\nu)|\mathscr{Y}^\nu_{[0,t]}\big)>0,
\ $ $Q^\nu$- and $Q^\beta$-a.s., so that
$$
\pi^{\beta\nu}_t(i)
=\frac{E\big(I(X^\nu_t=a_i)\frac{d\beta}{d\nu}(X^\nu_0)|
\mathscr{Y}^\nu_{[0,t]}\big)}{E\Big(\frac{d\beta}{d\nu}(X_0^\nu)|\mathscr{Y}^\nu_{[0,t]}\Big)}
$$
and it remains to notice that
$$
E\big(I(X^\nu_t=a_i)\frac{d\beta}{d\nu}(X^\nu_0)|
\mathscr{Y}^\nu_{[0,t]}\big)=\sum_{j=1}^n\frac{d\beta}{d\nu}(a_j)P\big(X^\nu_t=a_i,X^\nu_0=a_j|
\mathscr{Y}^\nu_{[0,t]}\big).
$$
Taking into  consideration \eqref{nuform}, we find
\begin{eqnarray*}
&&\big|\pi_t^\nu(i)-\pi^{\beta\nu}_t(i)\big|
=\Bigg|\pi_t^\nu(i)-\frac
{\sum_{j=1}^n\Big(\frac{d\beta}{d\nu}(a_j)P\big(X^\nu_0=a_j,
X_t^\nu=a_i|\mathscr{Y}^\nu_{[0,t]}\big)\Big)}{E\Big(\frac{d\beta}{d\nu}(X^\nu_0)|\mathscr{Y}^\nu_{[0,t]}
\Big)}\Bigg|
\\
&=&\frac{\Big|\sum_{j=1}^n\frac{d\beta}{d\nu}(a_j)\Big(
\pi^\nu_t(i)P\big(X^\nu_0=a_j|\mathscr{Y}^\nu_{[0,t]}\big)-
P\big(X^\nu_0=a_j,X^\nu_t=a_i|\mathscr{Y}^\nu_{[0,t]}\big)\Big)\Big|}
{E\Big(\frac{d\beta}{d\nu}(X^\nu_0)|\mathscr{Y}^\nu_{[0,t]}\Big)}.
\end{eqnarray*}
Then, since by the Jensen inequality $
1\big/E\big(\frac{d\beta}{d\nu}(X^\nu_0)|\mathscr{Y}^\nu_{[0,t]}\big)\le
E\big(\frac{d\nu}{d\beta}(X^\nu_0)|\mathscr{Y}^\nu_{[0,t]}\big), $
we get the chain of estimates
\begin{equation}\label{formula}
\begin{aligned}
& \big|\pi_t^\nu(i)-\pi^{\beta\nu}_t(i)\big|
\le\max_{a_j\in\mathbb{S}}\frac{d\beta}{d\nu}(a_j)
\max_{a_j\in\mathbb{S}}\frac{d\nu}{d\beta}(a_j)
\\
&\qquad\times
\Bigg|\sum_{j=1}^n\pi^\nu_t(i)\Big(P\big(X^\nu_0=a_j|\mathscr{Y}^\nu_{[0,t]}\big)-
P\big(X^\nu_0=a_j\big|X^\nu_t=a_i,\mathscr{Y}^\nu_{[0,t]}\big)\Big)\Bigg|
\\
&\quad\le\max_{a_j\in\mathbb{S}}\frac{d\beta}{d\nu}(a_j)
\max_{a_j\in\mathbb{S}}\frac{d\nu}{d\beta}(a_j)
\\
&\qquad\times\sum_{j=1}^n\pi^\nu_t(i)\Big|P\big(X^\nu_0=a_j|\mathscr{Y}^\nu_{[0,t]}\big)-
P\big(X^\nu_0=a_j|X^\nu_t=a_i,\mathscr{Y}^\nu_{[0,t]}\big)\Big|
\\
&\quad\le\max_{a_j\in\mathbb{S}}\frac{d\beta}{d\nu}(a_j)\max_{j
\in\mathbb{S}}\frac{d\nu}{d\beta}(a_j)
\\
&\qquad\times
\sum_{j=1}^n\Big|P\big(X^\nu_0=a_j|\mathscr{Y}^\nu_{[0,t]}\big)-
P\big(X^\nu_0=a_j|X^\nu_t=a_i,\mathscr{Y}^\nu_{[0,t]}\big)\Big|
\\
&\quad=\max_{a_j\in\mathbb{S}}\frac{d\beta}{d\nu}(a_j)\max_{a_j\in\mathbb{S}}\frac{d\nu}
{d\beta}(a_j)
\sum_{j=1}^n\Big|P\big(X^\nu_0=a_j|\mathscr{Y}^\nu_{[0,t]}\big)-\rho_{ji}(t)\Big|.
\end{aligned}
\end{equation}
The obvious formula $
P\big(X^\nu_0=a_j|\mathscr{Y}^\nu_{[0,t]}\big)=
\sum_{k=1}^n\pi^\nu_t(k)\rho_{jk}(t), $ and \eqref{formula}
provide
\begin{eqnarray}\label{new1}
\big|\pi_t^\nu(i)-\pi^{\beta\nu}_t(i)\big| &\le&
\max_{a_j\in\mathbb{S}}\frac{d\beta}{d\nu}(a_j)\max_{a_j\in\mathbb{S}}
\frac{d\nu}{d\beta}(a_j)
\sum_{j=1}^n\Big|\sum_{k=1}^n\pi^\nu_t(k)\rho_{jk}(t)-\rho_{ji}(t)\Big|
\nonumber\\
&\le&
\max_{a_j\in\mathbb{S}}\frac{d\beta}{d\nu}(a_j)\max_{a_j\in\mathbb{S}}
\frac{d\nu}{d\beta}(a_j)
\sum_{j=1}^n\sum_{k=1}^n\pi^\nu_t(k)\big|\rho_{jk}(t)-\rho_{ji}(t)\big|,
\end{eqnarray}
and \eqref{5.21m}.
Thus, by Lemma \ref{lem-3.1}, the desired statement \eqref{uspex2} holds true.

\subsection{The proof of Theorem \ref{theo-2.20}}
\label{sec-5.3}
We start with the following lemma.

\medskip
\begin{lemma}\label{lem-3.11}
Under the assumptions of Theorem {\rm \ref{theo-2.20}}, for any $t>0$
\begin{equation}\label{3.3a}
\max_{1\le j,k,\ell\le n}\big|\rho_{jk}(t)-\rho_{j\ell}(t)\big|\le
\exp\Big(-2t\min_{p\ne q}\sqrt{\lambda_{pq}\lambda_{qp}}\Big),
\end{equation}
\end{lemma}

\begin{proof}
Here we follow the notations from Lemma \ref{lem-3.1}. From
\eqref{bnd}, it follows that
\begin{eqnarray}\label{bnda}
\frac{d\triangle_t}{dt}\le
-\triangle_t\Bigg(\frac{\lambda_{i_\diamond i^\diamond}
\pi^\beta_t(i_\diamond)}{\pi^\beta_t(i^\diamond)} +
\frac{\lambda_{i^\diamond i_\diamond}
\pi^\beta_t(i^\diamond)}{\pi^\beta_t(i_\diamond)}\Bigg)
\end{eqnarray}
subject to $\triangle_0=1$. Set $\tau = \inf\{t:
i^\diamond(t)=i_\diamond(t)\}$. Since $\triangle_t$ is
nonincreasing function, $\triangle_t\equiv 0$ for  $t\ge \tau$,
and \eqref{3.3a} holds trivially. For $t<\tau$, as previously we find
$$
\begin{aligned}
\triangle_t &\le \exp\left\{-\int_0^t\Bigg(
\frac{\lambda_{i_\diamond i^\diamond}
\pi^\beta_s(i_\diamond)}{\pi^\beta_s(i^\diamond)} +
\frac{\lambda_{i^\diamond i_\diamond}
\pi^\beta_s(i^\diamond)}{\pi^\beta_s(i_\diamond)}\Bigg)ds
\right\}\\
&\le \exp\left\{-\int_0^t\min_{x\ge 0}\Bigg( \lambda_{i_\diamond
i^\diamond} x + \lambda_{i^\diamond i_\diamond}
\frac{1}{x}\Bigg)ds
\right\}\\
& = \exp\left\{-\int_0^t 2\sqrt{\lambda_{i_\diamond i^\diamond}
\lambda_{i^\diamond i_\diamond}} ds \right\}\le \exp\Big(
-2t\min_{p\ne q}\sqrt{\lambda_{pq}\lambda_{qp}} \Big),
\end{aligned}
$$
and \eqref{3.3a} follows.
\end{proof}

\bigskip
To prove the first statement of the theorem, taking into account $\nu\ll\beta$
we replicate a fragment from the
proof of Proposition \ref{pro-1.1}.

Using the notations introduced in section \ref{sec-2.1}, write
$\pi^\nu_t(i):=\pi^\nu_t(f)$ and
$\pi^{\beta\nu}_t(i):=\pi^{\beta\nu}_t(f)$ for $f(x)=I(x=a_i)$.
Then,
\begin{equation}\label{aliki}
E\big|\pi^{\beta\nu}_t(i)-\pi_t^\nu(i)\big| \le
E\Big|E\Big(\frac{d\nu}{d\beta}(X_0^\beta)\big|\mathscr{Y}^\beta_{[0,t]}\Big)
-
E\Big(\frac{d\nu}{d\beta}(X^\beta_0)\big|\mathscr{Y}^\beta_{[0,\infty)}
\vee \mathscr{X}^\beta_{[t,\infty)}\Big)\Big|
\end{equation}
and, since $(X^\beta,Y^\beta)$ is a Markov process,
$$
E\Big(\frac{d\nu}{d\beta}(X^\beta_0)\big|\mathscr{Y}^\beta_{[0,\infty)}
\vee \mathscr{X}^\beta_{[t,\infty)}\Big)=
E\Big(\frac{d\nu}{d\beta}(X^\beta_0)\big|\mathscr{Y}^\beta_{[0,t]}
\vee \mathscr{X}^\beta_t \Big).
$$
Then,
\begin{equation}\label{5.18m}
\begin{aligned}
&
E\Big(\frac{d\nu}{d\beta}(X_0^\beta)\big|\mathscr{Y}^\beta_{[0,t]}\Big)
-
E\Big(\frac{d\nu}{d\beta}(X^\beta_0)\big|\mathscr{Y}^\beta_{[0,\infty)}
\vee \mathscr{X}^\beta_{[t,\infty)}\Big)
\\
&= \sum_{j=1}^n\frac{d\nu}{d\beta}(a_j)\Big(P\big(X^\beta_0=a_j
|\mathscr{Y}^\beta_{[0,t]}\big) -P\big(X^\beta_0=a_j
|\mathscr{Y}^\beta_{[0,t]}\vee\mathscr{X}^\beta_t\big)\Big)
\\
&=\sum_{j=1}^n\sum_{\ell=1}^nI(X^\beta_t=a_\ell)\frac{d\nu}{d\beta}(a_j)
\Big(P\big(X^\beta_0=a_j|\mathscr{Y}^\beta_{[0,t]}\big)-\rho_{j\ell}(t)\Big)
\\
&=\sum_{j=1}^n\sum_{\ell=1}^n\sum_{k=1}^n\pi^\beta_t(k)
I(X^\beta_t=a_\ell)\frac{d\nu}{d\beta}(a_j)\big(\rho_{jk}(t)-\rho_{j\ell}(t)\big)
\\
&\le \max_{1\le j,k,\ell\le
n}\big|\rho_{jk}(t)-\rho_{j\ell}(t)\big|
\sum_{j=1}^n\frac{d\nu}{d\beta}(a_j).
\end{aligned}
\end{equation}

The first statement of Theorem \ref{theo-2.20} follows from
\eqref{aliki}, \eqref{5.18m}, and Lemma \ref{lem-3.11}.

\smallskip
The second statement follows from \eqref{5.21m} and Lemma \ref{lem-3.11}.

\section{Proofs for non-ergodic case} \label{sec-6}

Recall that in the non-ergodic setting under consideration
$$
\mathbb{S}=\Big\{\underbrace{a^1_1,\ldots,a^1_{n_1}}_{\mathbb{S}_1},
\ldots,\underbrace{a^m_1,\ldots,a^m_{n_m}}_{\mathbb{S}_m}\Big\}, \quad
m\ge 2
$$
with subalphabets $\mathbb{S}_1,\ldots,\mathbb{S}_m$
noncommunicating in the sense of \eqref{1.2'}.

\subsection{Auxiliary lemmas}
\label{sec-6.1}

In this subsection, $\widetilde{X}^j_t$ is an independent copy of
$X^j_t$ with the initial distribution $\mu^j$, defined on some
auxiliary probability space
$(\widetilde{\Omega},\widetilde{\mathscr{F}},\widetilde{P})$ and
$\widetilde{E}$ is the expectation with respect to
$\widetilde{P}$. Recall that $\mu^j$ is the invariant measure, so
that $\widetilde{X}^j_t$ is stationary process.

\begin{lemma} \label{lem-5.1}
Fix $r > 0$ and define
 $
 Z_n = \sum_{i=1}^n \big(Y^\beta_{ir} - Y^\beta_{(i-1)r}\big)^2.
 $
Then with $n \to \infty$
$$
 \frac{1}{n} Z_n \to r + \sum_{j=1}^m I(X^\beta_0
 \in \mathbb{S}_j)\widetilde{E}
 \left( \int_0^r h(\widetilde{X}^j_s)\,ds \right)^2.
$$
\end{lemma}

\begin{proof}
Define
 $$
 F(i) = E\Big[\Big( \int_0^r h(X^\beta_s)\,ds\Big)^2\Big|X^\beta_0=
 a_i\Big]
 $$
and $\mathscr{G}_n = \sigma\{Y_{[0,nr]}\} \vee
\sigma\{X_{[0,nr]}\}$.  Then
 $
 E\Big[\Big(Y^\beta_{(n+1)r} - Y^\beta_{nr}\Big)^2 \Big|
     \mathscr{G}_n \Big] = r + F(X^\beta_{nr})
 $
so that the sequence $M_n = Z_n - nr - \sum_{i=0}^{n-1}
F(X^\beta_{ir})$ is a martingale with respect to the filtration
$(\mathscr{G}_n)_{n\ge 1}$. It is easy to verify that there exists
$K < \infty$ such that for all $n$ we have $ E(M_{n+1}- M_n)^2 \le
K. $ It follows that $(1/n)M_n \to 0$ almost surely as $n \to
\infty$ (see, e.g., Chapter VII, Section 5, Theorem 4 in  \cite{Sh}).

Now consider $(1/n) \sum_{i=0}^{n-1} F(X^\beta_{ir})$.  If $X_0
\in \mathbb{S}_j$, then $X_t \in \mathbb{S}_j$ for all $t \ge 0$
and the process is ergodic in $\mathbb{S}_j$ with stationary
distribution $\mu^j$. Applying the ergodic theorem for each class
$\mathbb{S}_j$ we obtain
 $$
  \frac{1}{n} \sum_{i=0}^{n-1} F(X^\beta_{ir})  \to  \sum_{j=1}^m
  \widetilde{E}(F(\widetilde{X}_0))
           I(X_0 \in \mathbb{S}_j) = \sum_{j=1}^m
   \widetilde{E}\left(\int_0^r h(\widetilde{X}^j_s)\,ds\right)^2
            I(X^\beta_0 \in \mathbb{S}_j)
$$
as $n \to \infty$ a.s.
 Finally
  \begin{eqnarray*}
  \lim_{n \to \infty} \frac{1}{n} Z_n & =
     & \lim_{n\to \infty}\frac{1}{n} M_n
       + r +  \lim_{n\to\infty} \frac{1}{n} \sum_{i=0}^{n-1}
       F(X^\beta_{ir})\\
    & = & r + \sum_{j=1}^m \widetilde{E}
       \left(\int_0^r h(\widetilde{X}^j_s)\,ds \right)^2
           I(X^\beta_0 \in \mathbb{S}_j)
  \end{eqnarray*}
and we are done.
\end{proof}

With $\widetilde{X}^j_t$ defined as in Lemma \ref{lem-5.1} and
$r\ge 0$ let $ d_j(r) = \widetilde{E}\left(\int_0^r
h(\widetilde{X}^j_s)\, ds \right)^2. $

\begin{lemma} \label{lem-5.2} For any $k \neq j$ the following are equivalent:
\begin{enumerate}
\renewcommand{\theenumi}{{\rm \roman{enumi}}}
\item \label{i} $d_k(r) = d_j(r)$ for all $r \ge 0$; \item
\label{ii} $  h ^*_k\diag (\mu_k)\Lambda_k^q
 h _k =  h ^*_j\diag (\mu_j)\Lambda_j^q
 h _j $ for all $0 \le q \le n_i+n_j-1$.
\end{enumerate}
\end{lemma}
\begin{proof}
Notice first that
\begin{eqnarray*}
d_j(r) & = & 2\widetilde{E}\int_0^r \int_0^s h(\widetilde{X}^j_u)
h(\widetilde{X}^j_s)\,du\,ds
 = 2\int_0^r \int_0^s \widetilde{E}h(\widetilde{X}^j_u)
     h(\widetilde{X}^j_s)\,du\,ds
 \\
   & = & 2\int_0^r \int_0^s \widetilde{E}h(\widetilde{X}^j_0)
   h(\widetilde{X}^j_{s-u})\,du\,ds
= 2\int_0^r \int_0^s
    \widetilde{E}h(\widetilde{X}^j_0)h(\widetilde{X}^j_{v})\,dv\,ds.
 \end{eqnarray*}
Now, introduce  the vector $\widetilde{I}^j_t$ with entries $
I(\widetilde{X}^j_t=a_1^j),\ldots, I(\widetilde{X}^j_t=a_{n_j}^j)
$ and notice also that
\begin{eqnarray*}
\widetilde{E}h(\widetilde{X}^j_0)h(\widetilde{X}^j_{v})&=&
\widetilde{E} h ^*_j\widetilde{I}^j_0(\widetilde{I}^j_v)^* h _j
=\widetilde{E} h
^*_j\widetilde{I}^j_0(\widetilde{I}^j_0)^*e^{\Lambda_jv}
 h _j
\\
&=& h ^*_j\widetilde{E}\diag (\widetilde{I}^j_0)e^{\Lambda_jv}  h
_j = h ^*_j\diag (\mu^j)e^{\Lambda_jv} h _j.
\end{eqnarray*}
Therefore $ d_j(r)=2\int_0^r\int_0^s h ^*_j\diag
(\mu^j)e^{\Lambda_jv} h _j dvds $ so, $d_j(0) = {d_j}'(0) = 0$ and
$$
 {d_j}''(r) = 2 h ^*_j\diag (\mu^j)e^{\Lambda_jr} h _j.
$$
Differentiating with respect to $r$ a further $q$ times and then
putting $r=0$ we get
  $$
  d_j^{(2+q)}(0) =2 h ^*_j\diag (\mu^j)\Lambda^q_j h _j.
  $$
It follows immediately that if $d_k(r) = d_j(r)$ for all $r \ge
0$, then
$$
 h ^*_k\diag (\mu^k)\Lambda^q_k h _k
= h ^*_j\diag (\mu^j)\Lambda^q_j h _j
$$
for all $q \ge 0$ and so in particular for all $0 \le q \le
n_k+n_j-1$.

Suppose conversely that $  h ^*_j\diag (\mu^j)\Lambda^q_j h _j=  h
^*_k\diag (\mu^k)\Lambda^q_k h _k $ for all $0 \le q \le
n_k+n_j-1$. The Cayley$-$Hamilton theorem applied to the $(n_k+n_j)
\times (n_k +n_j)$ block diagonal matrix
   $
\left(\begin{smallmatrix}
   \Lambda_k & 0 \\ 0 & \Lambda_j
\end{smallmatrix}
\right)
   $
gives constants $c_0, c_1, \ldots , c_{n_k+n_j-1}$ so that
   $$
   \Lambda_k^{n_k+n_j} = \sum_{q = 0}^{n_k+n_j-1} c_q \Lambda_k^q
  \quad \mbox{and} \quad
   \Lambda_j^{n_i+n_j} = \sum_{q = 0}^{n_k+n_j-1} c_q \Lambda_j^q.
   $$
Therefore we have $ h ^*_k\diag(\mu^k)\Lambda_k^q
 h _k = h ^*_j\diag(\mu^j)\Lambda_j^q
 h _j$ for all $q>n_j+n_k-1$ as well. Using the fact that
$ e^{\Lambda_jr} = \sum_{q=0}^\infty \frac{r^q \Lambda_j^q}{q!}, $
we see that
 ${d_k}''(r) = {d_j}''(r)$ for all $r \ge 0$, and hence  $d_k(r) = d_j(r)$
 for all $r \ge 0$.
\end{proof}
\begin{lemma}\label{lem-4.3}
Assume {\rm(A-2)}. For any $\beta$
$$
\lim_{t\to\infty}E\Big|P\big(X^\beta_0\in\mathbb{S}_j|\mathscr{Y}^\beta_{[0,t]}\big)-
I(X^\beta_0\in\mathbb{S}_j)\Big|=0, \ j\ge 1.
$$
\end{lemma}

\begin{proof}
We use the notation $Z_n^{(r)}$ to express the dependence on $r$
of the function $Z_n$ in Lemma \ref{lem-5.1}.  We have $
\frac{1}{n}Y^\beta_n \to \sum_{j=1}^m h
^*_j\mu^jI(X^\beta_0\in\mathbb{S}_j) $ and
 $$
\frac{1}{n}Z_n^{(r)} \to r+\sum_{j=1}^md_j(r)I(X^\beta_0\in
\mathbb{S}_j)
$$ as $n \to \infty$, a.s.  Using the assumption (A-2) and
Lemma \ref{lem-5.2} we can find an integer $\ell$ and numbers
$r_i>0, i = 1,\ldots, \ell$ and construct a random variable of the
form $V_n = (Y^\beta_n,Z_n^{(r_1)}-nr_1,
\ldots,Z_n^{(r_\ell)}-nr_\ell)$ so that $ \frac{1}{n} V_n \to
\sum_{j=1}^m v_j I(X^\beta_0\in \mathbb{S}_j) $ as $n \to \infty$,
$P$-a.s, where the $v_1,\ldots,v_m$ are distinct vectors in
$\mathbb{R}^{\ell + 1}$.  Therefore $\{X^\beta_0\in
\mathbb{S}_j\}$ is $Y^\beta_{[0,\infty)}$-measurable a.s. and the
result follows immediately.
\end{proof}

\subsection{The proof of Theorem \ref{theo-2.2}}
\label{sec-6.2}

By Proposition \ref{pro-1.1}, it suffices to show that
 $$
\lim_{t\to\infty}E\big\|\pi^\beta_t-\pi^{\beta_0}_t\big\|=0.
 $$
We introduce a new filter, intermediate between $\pi^\beta_t$ and
$\pi_t^{\beta_0}$.   Define the random variable $U$ by $U = j$ on
the set $\{X_0^\beta \in \mathbb{S}_j\}$, and then define
  $$
  \pi^{\beta,U}_t(i) = P(X_t^\beta = a_i | \mathscr{Y}^\beta_{[0,t]},U).
  $$
Then
 \begin{eqnarray*}
\big\|\pi^\beta_t - \pi^{\beta,U}_t\big\| & = &
   \sum_{i=1}^n\left|P(X_t^\beta = a_i |\mathscr{Y}^\beta_{[0,t]})
      - P(X_t^\beta = a_i | \mathscr{Y}^\beta_{[0,t]},U)\right| \\
   & = &
   \sum_{i=1}^n\left|\sum_{j=1}^m P(X_t^\beta = a_i |\mathscr{Y}^\beta_{[0,t]},U=j)
   \left(P(U = j| \mathscr{Y}^\beta_{[0,t]})-I(U =
   j)\right) \right| \\
      & \le &
      \sum_{j=1}^m \left|P(U=j | \mathscr{Y}^\beta_{[0,t]})- I(U=j) \right|
\end{eqnarray*}
and
$$
\begin{aligned}
& \big\|\pi^{\beta,U}_t - \pi^{\beta_0}_t\big\|=
 \sum_{i=1}^n\left|P(X_t^\beta = a_i | \mathscr{Y}^\beta_{[0,t]},U) -
        P(X_t^\beta = a_i |\mathscr{Y}^\beta_{[0,t]},X_0^\beta) \right|
        \\
        &=\sum_{i=1}^n \sum_{j=1}^m I(U = j)
   \left|P(X_t^\beta = a_i | \mathscr{Y}^\beta_{[0,t]},U = j) -
        P(X_t^\beta = a_i |\mathscr{Y}^\beta_{[0,t]},U = j, X_0^\beta) \right|
        \\
        &=\sum_{j=1}^m I(U=j)\big\|\pi^{\beta^j}_t - \pi^{\beta^j_0}_t\big\|,
\end{aligned}
$$
where $\beta^j$ denotes the conditional distribution of $\beta$
restricted to the subalphabet $\mathbb{S}_j$. By Lemma
\ref{lem-4.3},
$$
\sum_{j=1}^m \left|P(U=j |
\mathscr{Y}^\beta_{[0,t]})- I(U=j) \right|\xrightarrow[t\to\infty]{P}0
$$
while
$\sum_{j=1}^m
I(U=j)\big\|\pi^{\beta^j}_t - \pi^{\beta^j_0}_t\big\|
\xrightarrow[t\to\infty]{\mathbb{L}_1}0$
by applying Theorem \ref{theo-2.0} to each
$\mathbb{S}_j$.

\appendix
\section{Proof of Proposition \ref{theo-A.2}}\label{app}

\begin{proof} {\em (Sketch)}
We use the following construction for $X$.
Let $X_0$ be a random variable with values in
$\mathbb{S}=\{1,2,3,4\}$ and $P(X_0=j)=\nu_j$, $j=1,\ldots,4$.
Introduce independent of $X_0$ matrix-valued process
\begin{equation}\label{4.16m}
\mathcal{N}_t=
  \begin{pmatrix}
    -N_{12}(t) & N_{12}(t) & 0 & 0\\
    0 & -N_{23}(t) & N_{23}(t) & 0\\
    0  & 0 & -N_{34}(t) & N_{34}(t) \\
    N_{41}(t) & 0 & 0 & -N_{41}(t)
  \end{pmatrix},
\end{equation}
where $N_{ij}(t)$ are independent copies of Poisson process with
the unit rate. Let us consider the It\^o equation
\begin{equation}\label{Ito.m}
I_t=I_0+\int_0^td\mathcal{N}^*_sI_{s-}
\end{equation}
with $I_0$ the vector with entries $I_0(j)=I(X_0=j)$,
$j=1,\ldots,4$. Since the jumps of Poisson processes $N_{ij}(t)$'s are
disjoint, for any $t>0$ the vector $I_t$ has only one nonzero
entry. Moreover, whereas the increments of $\mathcal{N}_t$ are
independent for nonoverlapping intervals, $I_t$ is Markov process.
It is  readily checked that, with the row vector $g=
  \begin{pmatrix}
    1 & 2 & 3 &4
  \end{pmatrix},
$ $X_t=g I_t$ is Markov process with values in $\mathbb{S}$ and
the transition intensities matrix $\Lambda$  and $I_t(j)=I(X_t=j)$,
$j=1,\ldots,4$.

We will follow Theorem 4.10.1  from \cite{LSMar}. The random
process $Y$ has piecewise constant paths with jumps of two
magnitudes, $+1$ and $-1$. Due to \eqref{Ito.m}, its saltus
measure $p(dt,dy)$ is completely described by
$$
\begin{aligned}
& p(dt,\{1\})=\big\{I_{t-}(4)dN_{41}(t)+I_{t-}(2)dN_{23}(t)\big\}
\\
&
p(dt,\{-1\})=\big\{I_{t-}(1)dN_{12}(t)+I_{t-}(3)dN_{34}(t)\big\}.
\end{aligned}
$$
So, the compensator $\overline{q}(dt,dy)$ of $p(dt,dy)$ with
respect to the filtration $(\mathscr{Y}_{[0,t]})_{t\ge 0}$ is
defined as
\begin{equation}\label{A.6m}
\begin{aligned}
\overline{q}(dt,\{1\})&=\big(\pi_{t-}(4)+\pi_{t-}(2)\big)dt=(1-Y_{t-})dt
\\
\overline{q}(dt,\{-1\})&=\big(\pi_{t-}(1)+\pi_{t-}(3)\big)dt=Y_{t-}dt.
\end{aligned}
\end{equation}
Notice also that
\begin{equation}\label{A.7m}
p(dt,\{1\})=(1-Y_{t-})dY_t\quad\text{and}\quad
p(dt,\{-1\})=-Y_{t-}dY_{t}.
\end{equation}
Equation \eqref{Ito.m} also gives ``drift+martingale''
presentation for $I_1(t)$, $I_2(t)$:
\begin{equation}\label{two}
\begin{aligned}
dI_t(1)&=\big(-I_{t}(1)+I_{t}(4)\big)dt+dM_1(t)
\\
dI_{t}(2)&=\big(I_{t}(1)-I_{t}(2)\big)dt+dM_2(t)
\end{aligned}
\end{equation}
with martingales
$$
\begin{aligned}
M_1(t)&=\int_0^t\Big(-I_{s-}(1)d(N_{12}(s)-s)+I_{s-}(4)d(N_{41}(s)-s)\Big)
\\
M_2(t)&=\int_0^t\Big(I_{s-}(1)d(N_{12}-s)-I_{s-}(2)d(N_{23}(s)-s)\Big).
\end{aligned}
$$
Then, by Theorem 4.10.1 in \cite{LSMar}, adapted to the case
considered, we have
\begin{equation}\label{A.9m}
\begin{aligned}
d\pi_1(t)&=\big(-\pi_t(1)+\pi_t(4)\big)dt +\int
H_1(\omega,t,y)\big[p(dt,dy)-\overline{q}(dt,dy)\big]
\\
d\pi_2(t)&=\big(\pi_t(1)-\pi_t(2)\big)dt +\int
H_2(\omega,t,y)\big[p(dt,dy)-\overline{q}(dt,dy)\big],
\end{aligned}
\end{equation}
where $H_i(\omega,t,y)$, $i=1,2$, are $ \mathscr{P}(Y)\otimes
\mathscr{B}(\mathbb{R}) $ -measurable functions (here
$\mathscr{B}(\mathbb{R})$ is the Borel $\sigma$-algebra on
$\mathbb{R}$ and  $\mathscr{P}(Y)$ is the predictable $\sigma $-algebra on
$\Omega \times \mathbb{R}_+$\ with respect to the filtration
$(\mathscr{Y}_{[0,t]}) _{t\ge 0}$). Moreover
$$
H_i(\omega,t,y)={\sf M}\big(\triangle M_i+I_-(i)|
\mathscr{P}(\mathcal{Y})\otimes
\mathscr{B}(\mathbb{R})\big)(\omega,t,y) -\pi_{t-}(i),
$$
where $\triangle M_i$ and $I_-(i)$ are the processes
$M_i(t)-M_i(t-)$ and $I_{t-}(i)$, respectively, and $ {\sf
M}\big(\cdot| \mathscr{P}(Y)\otimes \mathscr{B}(\mathbb{R})\big) $
is the conditional expectation with respect to the measure ${\sf
M}(d\omega,dt,dy)=P(d\omega)p(dt,dy)$ given $\mathscr{P}(Y)\otimes
\mathscr{B}(\mathbb{R}). $

By \eqref{two}, $\triangle M_i(t)+I_{t-}(i)=I_t(i)$ and the
structure of compensator $\overline{q}$ provides (here $\triangle
I_t(i)=I_t(i)- I_{t-}(i)$)
$$
{\sf
M}\big(I(i)|\mathscr{P}(Y)\otimes\mathscr{B}(\mathbb{R})\big)-\pi_{t-}(i)=
{\sf M}\big(\triangle
I(i)|\mathscr{P}(Y)\otimes\mathscr{B}(\mathbb{R})\big).
$$
The desired conditional expectation is determined uniquely from
the following identity: for any bounded, compactly supported in $t$ and
$\mathscr{P}(\mathcal{Y})\otimes
\mathscr{B}(\mathbb{R})$-measurable function $\phi(\omega,t,y)$
\begin{multline*}
E\int_0^\infty\int\phi(\omega,t,y)\triangle I_t(i)p(dt,dy)
\\
= E\int_0^\infty\int\phi(\omega,t,y) {\sf M}\big(\triangle
I(i)|\mathscr{P}(Y)\otimes \mathscr{B}(\mathbb{R})
\big)(\omega,t,y)\overline{q}(dt,dy).
\end{multline*}
By \eqref{Ito.m}
$$
\begin{aligned}
\triangle I_t(1)&=-I_{t-}(1)\triangle N_{12}(t)+I_{t-}(4)\triangle
N_{41}(t),
\\
\triangle I_t(2)&=I_{t-}(1)\triangle N_{12}(t)-I_{t-}(2)\triangle
N_{23}(t),
\end{aligned}
$$
and so
$$
\begin{aligned}
& \triangle I_t(1)p(dt,\{1\})=I_{t-}(4)dN_{41}(t),
\\
& \triangle I_t(1)p(dt,\{-1\})=-I_{t-}(1)dN_{12}(t),
\\
& \triangle I_t(2)p(dt,\{1\})=-I_{t-}(2)dN_{23}(t),
\\
& \triangle I_t(2)p(dt,\{-1\})=I_{t-}(1)dN_{12}(t).
\end{aligned}
$$
Owing to the obvious relations
\begin{gather*}
I_4(t)\equiv I_4(t)(1-Y_t), \ \ I_2(t)\equiv I_2(t)(1-Y_t),
\\
I_1(t)\equiv I_1(t)Y_t, \ \ I_3(t)\equiv I_3(t)Y_t
\end{gather*}
we have
\begin{equation}\label{A.10m}
\begin{aligned}
\pi_{t-}(2)dt=\pi_{t-}(2)(1-Y_{t-})dt, \ \
\pi_{t-}(2)dt=\pi_{t-}(2)(1-Y_{t-})dt
\\
\pi_{t-}(1)dt=\pi_{t-}(1)Y_{t-}dt, \ \
\pi_{t-}(3)dt=\pi_{t-}(3)Y_{t-}dt.
\end{aligned}
\end{equation}
Taking into account \eqref{A.6m}, we find
$$
\begin{aligned}
& H_1(\omega,t,y)=
  \begin{cases}
    \pi_{t-}(4),& y=1, \\
    -\pi_{t-}(1), & y=-1,
  \end{cases}
\\
& H_2(\omega,t,y)=
  \begin{cases}
    -\pi_{t-}(2),& y=1, \\
    \pi_{t-}(1), & y=-1.
  \end{cases}
\end{aligned}
$$
In accordance with \eqref{A.6m}, \eqref{A.7m}, the formulae for $H_1, \
H_2$, and \eqref{A.10m}, we transform \eqref{A.9m} to
$$
\begin{aligned}
d\pi_1(t)&=\big(-\pi_t(1)+\pi_t(4)\big)dt
+\pi_{t-}(4)(1-Y_{t-})(dY_t-dt)+\pi_{t-}(1)Y_{t-}(dY_t+dt)
\\
&=\pi_{t-}(4)(1-Y_{t-})dY_t+\pi_{t-}(1)Y_{t-}dY_t
\\
&=\big(1-\pi_{t-}(2)\big)(1-Y_{t-})dY_t+\pi_{t-}(1)Y_{t-}dY_t,
\\
d\pi_2(t)&=\big(\pi_t(1)-\pi_t(2)\big)dt
-\pi_{t-}(2)(1-Y_{t-})(dY_t-dt)-\pi_{t-}(1)Y_{t-}(dY_t+dt)
\\
&=-\pi_{t-}(2)(1-Y_{t-})dY_t-\pi_{t-}(1)Y_{t-}dY_t.
\end{aligned}
$$
\end{proof}

\bigskip

{\bf Acknowledgements.} The authors gratefully acknowledge Boris Tsirelson
for bringing \cite{Wei} and  the example in \cite{Wil} to their
attention,  Rami Atar for suggesting use of Theorem 1 in
\cite{AZ2} for the proof of Theorem \ref{theo-2.0}, and the anonymous
referees whose comments and advises allowed us to improve
the paper significantly.
%----------------------------------------------------------------

\end{document}